\documentclass[letterpaper, 10 pt, conference]{ieeeconf}  

\IEEEoverridecommandlockouts                              

\usepackage{amsmath}
\usepackage{amssymb}
\usepackage{tikz}
\usetikzlibrary{automata,positioning}
\usepackage{pgfplots}
\usepackage{graphicx}
\usepackage{setspace}
\usepackage{caption}
\usepackage{subcaption}
\usepackage{xcolor}
\usepackage{mathtools}
\usepackage{epstopdf}

\newcommand{\until}[2]{\, U_{[{#1},{#2}]} \, }

\newcommand{\rsxt}[2]{\rho^{#1}(\boldsymbol{x}_\phi,{#2})} 
\newcommand{\rsct}[2]{\rho^{#1}(\boldsymbol{x}_{#1},{#2})} 
\newcommand{\rsnt}[1]{\rho^{#1}(\boldsymbol{x}_{\phi_i})}
\newcommand{\rsntk}[1]{\rho^{#1}(\boldsymbol{x}_{\phi_{\mathfrak{c}_i}})}
\newcommand{\rsntnc}[1]{\rho^{#1}(\boldsymbol{x}_i)}
\newcommand{\rso}[1]{\rho^{\text{opt}}_{i}}
\newcommand{\rsntont}[1]{\rho^{#1}(\boldsymbol{x}_{\phi_i}(0))}
\newcommand{\rsnttnt}[1]{\rho^{#1}(\boldsymbol{x}_{\phi_i}(t))}

\newcommand{\rsnttm}[1]{\rho^{#1}\big(\boldsymbol{x}_{\phi_i}(t)\big)}
\newcommand{\rsntj}[1]{\rho^{#1}(\boldsymbol{x}_{\phi_j})}

\newtheorem{definition}{Definition} 
\newtheorem{theorem}{Theorem} 
\newtheorem{assumption}{Assumption} 
\newtheorem{problem}{Problem} 
\newtheorem{remark}{Remark}
\newtheorem{lemma}{Lemma}
\newtheorem{corollary}{Corollary}
\newtheorem{example}{Example}

\title{\LARGE \bf
Decentralized Robust Control of Coupled Multi-Agent Systems \\under Local Signal Temporal Logic Tasks
}

\author{Lars Lindemann and Dimos V. Dimarogonas
\thanks{This work was supported in part by the Swedish Research Council (VR), the European Research Council (ERC), the Swedish Foundation for Strategic Research (SSF), the EU H2020 Co4Robots project, the SRA ICT TNG project STaRT, and the Knut and Alice Wallenberg Foundation (KAW).}
\thanks{The authors are with the Department of Automatic Control, School of Electrical Engineering, KTH Royal Institute of Technology, 100 44 Stockholm, Sweden. {\tt\small llindem@kth.se (L. Lindemann), dimos@kth.se (D.V. Dimarogonas)}}%
}

\begin{document}

\maketitle
\thispagestyle{empty}
\pagestyle{empty}

\begin{abstract}

Motivated by the recent interest in formal methods-based control of multi-agent systems, we adopt a bottom-up approach. Each agent is subject to a local signal temporal logic task that may depend on other agents' behavior. These dependencies pose control challenges since some of the tasks may be opposed to each other. We first develop a local continuous feedback control law and identify conditions under which this control law guarantees satisfaction of the local tasks. If these conditions do not hold, we propose to use the developed control law in combination with an online detection \& repair scheme, expressed as a local hybrid system. After detection of a critical event, a three-stage procedure is initiated to resolve the problem. The theoretical results are illustrated in simulations.

\end{abstract}

\section{Introduction}
\label{sec:introduction}

Multi-agent systems under global objectives such as consensus, formation control, and connectivity maintenance have been well studied by the research community. Comprehensive overviews of these topics can be found in \cite{cao2013overview} and \cite{mesbahi2010graph}, where the derived controllers are mainly distributed control laws. The need for more complex and rich objectives in robotic applications has led to formal methods-based control strategies where temporal logics, e.g., linear temporal logic, are used to formulate high-level temporal tasks. Top-down approaches have been considered in \cite{kloetzer2010automatic,nikou1} by decomposing a global temporal task into local ones that need to be executed by each agent individually. Top-down approaches often require agent synchronization and are usually subject to high computational complexity and hence impractical when the problem size becomes larger. On the other hand, the works in \cite{guo2013reconfiguration,tumova2016multi} favor a bottom-up approach, where local tasks are independently distributed to each agent. This leads to partially  decentralized solutions that reduce the computational burden. In a bottom-up approach, feasibility of each local task does not imply feasibility of the conjunction of all local tasks \cite{guo2013reconfiguration} since some of the local tasks may be opposed to each other. The presented works in \cite{kloetzer2010automatic}-\cite{tumova2016multi} rely on automata-based verification techniques that discretize the physical environment and agent dynamics. In this paper, we instead consider continuous-time and nonlinear dynamics without the need for discretizing neither environment nor agent dynamics in space or time. To the best of our knowledge, this is the first approach not making use of such discretization in the context of formal methods-based multi-agent control. This paper extends our work on single-agent systems \cite{lindemann2017prescribed} to multi-agent systems.

We adopt a bottom-up approach by considering local tasks formulated in signal temporal logic \cite{maler1}. These tasks can depend on each other, i.e., also oppose each other. This makes the control of multi-agent systems under signal temporal logic tasks a challenge and the main research question in this paper. Signal temporal logic introduces the notion of space robustness \cite{donze2}, a robustness metric stating how robustly a signal satisfies a given task. In a first step, we identify conditions under which a continuous feedback control law, which is derived by combining space robustness and prescribed performance control \cite{bechlioulis2014low}, satisfies basic signal temporal logic tasks. If these conditions do not hold, an online detection \& repair scheme is introduced by defining a local hybrid system \cite{goebel2012hybrid} for each agent. Critical events will be detected and resolved in a three-stage procedure, gradually relaxing parameters such as robustness. One advantage of our decentralized approach is the low computational complexity due to the continuous feedback control laws. Furthermore, the team of agents is allowed to be heterogeneous with additional dynamic couplings among them. Robustness is considered with respect to disturbances and with respect to the signal temporal logic task. Multi-agent systems under signal temporal logic tasks have also been considered in \cite{liu2017communication} in a centralized approach, not investigating formula dependencies, but with a special focus on communication. 

The remainder is organized as follows: in Section \ref{sec:preliminaries}, notation and preliminaries are introduced, while Section \ref{sec:problem_statement} presents the problem definition. Section \ref{sec:inf} presents our solution to the stated problem, which is verified by simulations in Section \ref{sec:sim}. Conclusions are given in Section \ref{sec:conclusion}. This online version is an extended version of the 2018 American Control Conference version.

\section{Preliminaries}
\label{sec:preliminaries}

Scalar quantities are denoted by lowercase, non-bold letters $x$ and column vectors are lowercase, bold letters $\boldsymbol{x}$. True and false are denoted by $\top$ and $\bot$; $\mathbb{R}$ are the real numbers, while $\mathbb{R}^n$ is the $n$-dimensional real vector space. The natural, non-negative, and positive real numbers are $\mathbb{N}$, $\mathbb{R}_{\ge0}$, and $\mathbb{R}_{>0}$, respectively. For convenience, we define $\begin{bmatrix}
\boldsymbol{x} & \boldsymbol{y}
\end{bmatrix}:=\begin{bmatrix}
\boldsymbol{x}^T & \boldsymbol{y}^T
\end{bmatrix}^T$. For two sets $\mathcal{X}$ and $\mathcal{Y}$, the set-valued map $F:\mathcal{X}\rightrightarrows\mathcal{Y}$ maps each $\boldsymbol{x}\in \mathcal{X}$ to a set $F(\boldsymbol{x})\subseteq\mathcal{Y}$. The inverse image by a function $F$ of a set $\mathcal{M}\subseteq\mathcal{Y}$  is given by inv$\big(F(\mathcal{M})\big):=\{\boldsymbol{x}\in\mathcal{X}|F(\boldsymbol{x})\cap \mathcal{M}\neq \emptyset\}$.

Two basic results regarding the existence of solutions for initial-value problems (IVP) are needed in this paper. Assume $\boldsymbol{y}\in\Omega_{\boldsymbol{y}}\subseteq \mathbb{R}^{n_{\boldsymbol{y}}}$ and consider the IVP
\begin{align}
\dot{\boldsymbol{y}} = H(\boldsymbol{y},t) \label{eq:xi_systems} \text{ with } \boldsymbol{y}_0:=\boldsymbol{y}(0)\in \Omega_{\boldsymbol{y}},
\end{align}
where $H:\Omega_{\boldsymbol{y}}\times \mathbb{R}_{\ge 0} \to \mathbb{R}^{n_{\boldsymbol{y}}}$ and $\Omega_{\boldsymbol{y}}$ is a non-empty and open set. A solution to this IVP is a signal $\boldsymbol{y}:\mathcal{J}\to\Omega_{\boldsymbol{y}}$ with $\mathcal{J}\subseteq \mathbb{R}_{\ge 0}$ obeying \eqref{eq:xi_systems}.

\begin{lemma}\cite[Theorem 54]{sontag2013mathematical}\label{theorem:sontag1}
Consider the IVP in \eqref{eq:xi_systems}. Assume that $H:\Omega_{\boldsymbol{y}}\times \mathbb{R}_{\ge 0}\to \mathbb{R}^{n_{\boldsymbol{y}}}$ is: 1) locally Lipschitz continuous on $\boldsymbol{y}$ for each $t\in \mathbb{R}_{\ge 0}$, 2) piecewise continuous on $t$ for each fixed $\boldsymbol{y}\in\Omega_{\boldsymbol{y}}$.  Then, there exists a unique and maximal solution $\boldsymbol{y}:\mathcal{J}\to\Omega_{\boldsymbol{y}}$ with $\mathcal{J}:=[0,\tau_{\text{max}})\subseteq \mathbb{R}_{\ge 0}$ and $\tau_{\text{max}}\in \mathbb{R}_{>0}\cup\infty$.
\end{lemma}

\begin{lemma}\cite[Proposition C.3.6]{sontag2013mathematical}\label{theorem:sontag2}
Assume that the assumptions of Lemma \ref{theorem:sontag1} hold. For a maximal solution $\boldsymbol{y}$ on $\mathcal{J}=[0,\tau_{\text{max}})$ with $\tau_{\text{max}}<\infty$ and for any compact set $\Omega_{\boldsymbol{y}}^\prime\subset\Omega_{\boldsymbol{y}}$, there exists $t^\prime\in \mathcal{J}$ such that $\boldsymbol{y}(t^\prime)\notin\Omega_{\boldsymbol{y}}^\prime$.
\end{lemma}  

\subsection{Signal Temporal Logic (STL)}
Signal temporal logic (STL) is a predicate logic based on signals \cite{maler1}.  STL consists of predicates $\mu$ that are obtained after the evaluation of a predicate function $h:\mathbb{R}^n\to\mathbb{R}$ as
\begin{align*}
	\mu:=
	\begin{cases} 
		\top \text{ if } h(\boldsymbol{x})\ge0\\
		\bot \text{ if } h(\boldsymbol{x})< 0.
	\end{cases}
\end{align*} 

For instance, it is possible to express the predicate $\mu:=(|{x}_{i}+{x}_{j}|\le 1)$ with the predicate function $h(\boldsymbol{x}):=1-|{x}_{i}+{x}_{j}|$ to specify that the $i$-th and $j$-th state should be close. The STL syntax is
\begin{align*}
	\phi \; ::= \; \top \; | \; \mu \; | \; \neg \phi \; | \; \phi \wedge \psi \; | \; \phi  \until{a}{b} \psi\;,
\end{align*}
where $\mu$ is a predicate and $\phi$ and $\psi$ are STL formulas. The temporal until-operator $\until{a}{b}$ is time bounded with time interval $[a,b]$ where $a,b\in \mathbb{R}_{\ge 0}$ is such that $a\le b$.  The satisfaction relation $(\boldsymbol{x},t) \models \phi$ indicates if the signal $\boldsymbol{x}:\mathbb{R}_{\ge 0}\to \mathbb{R}^n$ satisfies $\phi$ at time $t$. The STL semantics are given next.
\begin{definition}[STL Semantics]
The STL semantics  are inductively defined as \cite[Definition 1]{maler1}:
\begin{align*}
&(\boldsymbol{x},t) \models \mu \hspace{-0.4cm}&  &\Leftrightarrow h(\boldsymbol{x}(t))\ge 0 \\
&(\boldsymbol{x},t) \models \neg\phi\hspace{-0.4cm}&  &\Leftrightarrow \neg((\boldsymbol{x},t) \models \phi)\\
&(\boldsymbol{x},t) \models \phi \wedge \psi\hspace{-0.4cm}&  &\Leftrightarrow (\boldsymbol{x},t) \models \phi \wedge (\boldsymbol{x},t) \models \psi\\
&(\boldsymbol{x},t) \models \phi \until{a}{b} \psi\hspace{-0.4cm}&  &\Leftrightarrow	\exists t_1 \in[t+a,t+b]\text{ s.t. } (\boldsymbol{x},t_1)\models \psi \\ 
& \hspace{-0.4cm}& & \hspace{0.43cm}\wedge  \forall t_2\in[t,t_1],\; (\boldsymbol{x},t_2)\models \phi
\end{align*} 
\end{definition}

 Disjunction-, eventually-, and always-operator are derived as $\phi\vee\psi:=\neg(\neg\phi\wedge\neg\psi)$, $F_{[a,b]}\phi:=\top  \until{a}{b} \phi$, and $G_{[a,b]}\phi := \neg F_{[a,b]}\neg\phi$, respectively. Robust semantics, called space robustness and denoted by $\rho^{\phi}(\boldsymbol{x},t)$, have been introduced in \cite{donze2} and are defined in Definition \ref{def:2}; $\rho^{\phi}(\boldsymbol{x},t)$ determines how robustly the signal $\boldsymbol{x}:\mathbb{R}_{\ge 0}\to \mathbb{R}^n$ satisfies $\phi$ at time $t$. It holds that $(\boldsymbol{x},t)\models \phi$ if $\rho^{\phi}(\boldsymbol{x},t)>0$.
\begin{definition}[Space Robustness] {The semantics of space robustness are inductively defined as \cite[Definition 3]{donze2}:}
	\begin{align*}
		\rho^{\mu}(\boldsymbol{x},t)&:= h(\boldsymbol{x}(t))\\
		\rho^{\neg\phi}(\boldsymbol{x},t) &:= 	-\rho^{\phi}(\boldsymbol{x},t)\\
		\rho^{\phi \wedge \psi}(\boldsymbol{x},t) &:= 	\min\big(\rho^{\phi}(\boldsymbol{x},t),\rho^{\psi}(\boldsymbol{x},t)\big)\\
		\rho^{F_{[a,b]} \phi}(\boldsymbol{x},t) &:= \underset{t_1\in[t+a,t+b]}{\max}\rho^{\phi}(\boldsymbol{x},t_1)\\
		\rho^{G_{[a,b]}\phi}(\boldsymbol{x},t) &:= \underset{t_1\in[t+a,t+b]}{\min}\rho^{\phi}(\boldsymbol{x},t_1)
	\end{align*}
	\label{def:2}
\end{definition}


The definitions of $\rho^{\phi \vee \psi}(\boldsymbol{x},t)$ and $\rho^{\phi \until{a}{b} \psi}$ are omitted since they will not be considered in the remainder. We abuse the notation as $ \rho^\phi(\boldsymbol{x}(t)) :=\rho^{\phi}(\boldsymbol{x},t)$ if $t$ is not explicitly contained in $\rho^{\phi}(\boldsymbol{x},t)$. For instance, $\rho^\mu(\boldsymbol{x}(t))  :=\rho^{\mu}(\boldsymbol{x},t) := h(\boldsymbol{x}(t))$ since $h(\boldsymbol{x}(t))$ does not contain $t$ as an explicit parameter. However, $t$ is explicitly contained in $\rho^{\phi}(\boldsymbol{x},t)$ if and only if temporal operators (eventually, always, or until) are used.

\subsection{A Bottom-up Approach for Multi-Agent Systems}

Consider a multi-agent system that consists of $M$ agents and where each agent is possibly affecting the behavior of another agent. Therefore, communication among agents is crucial. We model the communication by using a static, i.e., time-independent, and undirected graph $\mathcal{G}:=(\mathcal{V},\mathcal{E})$ \cite{mesbahi2010graph}. The vertex set is $\mathcal{V}:=\{v_{1},v_{2},\hdots,v_{M}\}$, while the edge set is $\mathcal{E}\in \mathcal{V}\times \mathcal{V}$.  Two agents $v_{i},v_{j}\in\mathcal{V}$ can communicate if and only if there exists a path between $v_i$ and $v_j$. A path is a sequence  $v_i,v_{k_1},\hdots,v_{k_P},v_j$ such that $(v_i,v_{k_1}),\hdots,(v_{k_P},v_j)\in\mathcal{E}$. As a consequence, all agents can communicate if and only if $\mathcal{G}$ is connected.


Let $\boldsymbol{x}_{i}\in\mathbb{R}^n$, $\boldsymbol{u}_{i}\in\mathbb{R}^{m_i}$, and $\boldsymbol{w}_i\in \mathcal{W}_i$ be the state, input, and additive noise of agent $v_{i}$'s dynamics with $\mathcal{W}_i\subset\mathbb{R}^n$ being a bounded set. Let $\boldsymbol{x}:=\begin{bmatrix}{\boldsymbol{x}_{1}} & {\boldsymbol{x}_{2}} & \hdots & {\boldsymbol{x}_{M}}\end{bmatrix}$ be the stacked vector of all agents' states. Each agent $v_{i}$ obeys the nonlinear and coupled dynamics
\begin{align}\label{eq:system}
\dot{\boldsymbol{x}}_{i}=f_i(\boldsymbol{x}_i)+f_i^\text{c}(\boldsymbol{x})+g_i(\boldsymbol{x}_{i})\boldsymbol{u}_{i}+\boldsymbol{w}_i,
\end{align}
where $f_i^\text{c}(\boldsymbol{x})$ is a term describing preassumed dynamic couplings of the multi-agent system. Also define $\boldsymbol{x}_i^{\text{ext}}:=\begin{bmatrix} {\boldsymbol{x}_{j_1}} & \hdots & {\boldsymbol{x}_{j_{M-1}}}  \end{bmatrix}$ such that $v_{j_1},\hdots,v_{j_{M-1}}\in\mathcal{V}\setminus \{v_i\}$, i.e., $\boldsymbol{x}_i^{\text{ext}}$ is a stacked vector containing the states of all agents except of $\boldsymbol{x}_i$. Note that $\boldsymbol{x}_i^{\text{ext}}$ is contained in $f_i^\text{c}(\boldsymbol{x})$ and can be seen as an external input generated by an exo-system, i.e., other agents. The functions $f_i$, $f_i^\text{c}$, and $g_i$ need to satisfy Assumption \ref{assumption:1}.
\begin{assumption}\label{assumption:1}
The functions $f_i:\mathbb{R}^n\to\mathbb{R}^n$, $f_i^\text{c}:\mathbb{R}^{nM}\to\mathbb{R}^n$, and $g_{i}:\mathbb{R}^n\to\mathbb{R}^{n\times m_i}$ are locally Lipschitz continuous, and $g_i(\boldsymbol{x}_{i}){g_i(\boldsymbol{x}_{i}}^T)$ is positive definite for all $\boldsymbol{x}_{i}\in \mathbb{R}^n$.
\end{assumption}  
\begin{remark}
The term $f_i^\text{c}(\boldsymbol{x})$ represents preassumed dynamic couplings that the multi-agent system is subject to. These couplings can, for instance, express consensus, formation control, connectivity maintenance, or obstacle avoidance objectives. 
\end{remark}

We now tailor the definitions of STL and its robust semantics to multi-agent systems. In our bottom-up approach, each agent $v_{i}\in\mathcal{V}$ is subject to a local STL formula. As a notational rule, the local formula of agent $v_i$ is endowed with the subscript, i.e., $\phi_i$. Based on \cite[Definition 3]{tumova2016multi}, local satisfaction of $\phi_i$ by the signal $\boldsymbol{x}_{\phi_i}:\mathbb{R}_{\ge 0}\to\mathbb{R}^{p_i}$ is defined in Definition \ref{def:local_sat}. We will be more specific regarding $\boldsymbol{x}_{\phi_i}$ and $p_i$ after Definition \ref{def:dependency}. 
\begin{definition} [Local Satisfaction]\label{def:local_sat}
The signal $\boldsymbol{x}_{\phi_i}:\mathbb{R}_{\ge 0}\\ \to\mathbb{R}^{p_i}$ locally satisfies $\phi_i$ if and only if $(\boldsymbol{x}_{\phi_i},0)\models \phi_i$.
\end{definition} 

Local  feasibility of $\phi_i$ is next defined in Definition \ref{def:local_inf}.
\begin{definition} [Local Feasibility]\label{def:local_inf}
The formula $\phi_i$ is locally feasible if and only if $\exists \boldsymbol{x}_{\phi_i}:\mathbb{R}_{\ge 0}\to\mathbb{R}^{p_i}$ such that $\boldsymbol{x}_{\phi_i}$ locally satisfies $\phi_i$.
\end{definition} 

Each local formula $\phi_i$ depends on agent $v_i$ and may also depend on some other agents $v_j\in\mathcal{V}$. Consider $\boldsymbol{x}_j:\mathbb{R}_{\ge 0}\to\mathbb{R}^n$ to be the solution to \eqref{eq:system} associated with agent $v_j$.
\begin{definition}[Formula-Agent Dependency]\label{def:dependency}
If $\boldsymbol{x}_j(t)$ is not contained in $\boldsymbol{x}_{\phi_i}(t)$ for all $t\in\mathbb{R}_{\ge 0}$ and local satisfaction of $\phi_i$, i.e., $(\boldsymbol{x}_{\phi_i},0)\models \phi_i$, can be evaluated, then $\phi_i$ does not depend on $v_j$. Otherwise, i.e., knowledge of $\boldsymbol{x}_j(t)$ is needed and hence $\boldsymbol{x}_j(t)$ is contained in $\boldsymbol{x}_{\phi_i}(t)$, then $\phi_i$ does depend on $v_j$ and we say that agent $v_j$ is participating in $\phi_i$. 
\end{definition}

The set of participating agents in $\phi_i$ is 
\begin{align*}
\mathcal{V}_{\phi_i}:=\{v_{j_1},v_{j_2},\hdots,v_{j_{P(\phi_i)}}\}\subseteq\mathcal{V},
\end{align*} 
where $P(\phi_i):=\sum_{j=1}^{|\mathcal{V}|} \chi_j(\phi_i)$ is a function evaluating the total number of participating agents in $\phi_i$ with 
\begin{align*}
\chi_j(\phi_i):=
 \begin{cases} 
 1 &\text{if } \phi_i \text{ depends on } v_j \\
 0 &\text{otherwise.}
 \end{cases}
\end{align*} 
It holds that each $v_j\in\mathcal{V}_{\phi_i}$ is participating in $\phi_i$ and $\nexists v_k\in \mathcal{V}\setminus\mathcal{V}_{\phi_i}$ such that $v_k$ is participating in $\phi_i$. Define 
\begin{align*}
\boldsymbol{x}_{\phi_i}(t):=\begin{bmatrix} {\boldsymbol{x}_{j_1}(t)} & \hdots & {\boldsymbol{x}_{j_{P(\phi_i)}}(t)} \end{bmatrix} \end{align*} 
for all $t\in\mathbb{R}_{\ge 0}$ with $v_{j_1},\hdots,v_{j_{P(\phi_i)}}\in \mathcal{V}_{\phi_i}$, i.e., all agents participating in $\phi_i$. Finally, for the signal $\boldsymbol{x}_{\phi_i}:\mathbb{R}_{\ge 0}\to\mathbb{R}^{p_i}$ we conclude that $p_i:=nP(\phi_i)$ so that $\boldsymbol{x}_{\phi_i}$ is completely defined. 

We call $\phi_i$ a non-collaborative formula if and only if $P(\phi_i)=1$. In other words, the satisfaction of $\phi_i$ does not depend on other agents $v_{j}\in\mathcal{V}\setminus \{v_i\}$, and hence $\boldsymbol{x}_{\phi_i}=\boldsymbol{x}_i$. Otherwise, i.e., if $P(\phi_i)> 1$, we call $\phi_i$ a collaborative formula. Since $\phi_i$ always depends on $v_i$, it always holds that $P(\phi_i)\ge 1$. Global satisfaction of the set of formulas $\{\phi_1,\hdots,\phi_M\}$ by the signal $\boldsymbol{x}:\mathbb{R}_{\ge 0}\to\mathbb{R}^{nM}$ is introduced in Definition \ref{def:glob_sat}. Note that $\boldsymbol{x}_{\phi_i}$ is naturally contained in $\boldsymbol{x}$.
\begin{definition}[Global Satisfaction]\label{def:glob_sat}
The signal $\boldsymbol{x}:\mathbb{R}_{\ge 0 }\to \mathbb{R}^{nM}$ globally satisfies $\{\phi_1,\hdots,\phi_M\}$ if and only if $\boldsymbol{x}_{\phi_i}$ locally satisfies $\phi_i$ for all agents $v_i\in\mathcal{V}$.
\end{definition}

In this respect, we similarly define global feasibility.
\begin{definition}[Global Feasibility]
The set of formulas $\{\phi_1,\hdots,\phi_M\}$ is globally feasible if and only if $\exists \boldsymbol{x}:\mathbb{R}_{\ge 0}\to\mathbb{R}^{nM}$ such that $\boldsymbol{x}$ globally satisfies $\{\phi_1,\hdots,\phi_M\}$.
\end{definition}

Next, maximal dependency clusters are introduced in a similar vein as in \cite[Definition 4]{guo2013reconfiguration}. 
\begin{definition}[Maximal Dependency Cluster]\label{def:max_cluster}
Consider the undirected dependency graph $\mathcal{G}_d:=(\mathcal{V},\mathcal{E}_d)$ where there is an edge $(v_i,v_j) \in \mathcal{E}_d\subseteq \mathcal{V}\times \mathcal{V}$ if and only if the formula $\phi_i$ depends on $v_j$ in the sense of Definition \ref{def:dependency}; $\Xi\subseteq \mathcal{V}$ is a maximal dependency cluster if and only if $\forall v_i,v_j\in \Xi$ there is a path from $v_i$ to $v_j$ in $\mathcal{G}_d$ and $\nexists v_i\in \Xi, v_k\in \mathcal{V}\setminus \Xi$ such that there is a path from $v_i$ to $v_k$. 
\end{definition}

Consequently, a multi-agent system under $\{\phi_1,\hdots,\phi_M\}$ induces $L\le M$ maximal dependency clusters denoted by $\bar{\Xi}:=\{\Xi_1,\hdots,\Xi_L\}$. These clusters are maximal in the sense that there are no formula-agent dependencies between clusters, i.e., $\nexists v_i\in\Xi_{l_1},v_j\in\Xi_{l_2}$ with $l_1\neq l_2$ such that $\phi_i$ depends on $v_j$, which is different to \cite[Definition 4]{guo2013reconfiguration}. Even though maximal dependency clusters have no formula-agent dependencies, dynamic couplings between clusters induced by $f_i^\text{c}(\boldsymbol{x})$  may be present.
\begin{example}
Consider three agents $v_1$, $v_2$, and $v_3$ with $\phi_1:=F_{[a_1,b_1]}(\|\boldsymbol{x}_1-\boldsymbol{x}_2\|\le 1)$, $\phi_2:=F_{[a_2,b_2]}(\|\boldsymbol{x}_2\|\le 1)$, and $\phi_3:=F_{[a_3,b_3]}(\|\boldsymbol{x}_3\|\le 1)$. Then $\bar{\Xi}:=\{\Xi_1,\Xi_2\}$ with $\Xi_1=\{v_1,v_2\}$ and $\Xi_2=\{v_3\}$. 
\end{example}

\subsection{Hybrid Systems}

Hybrid systems have recently been modeled and analyzed in \cite{goebel2012hybrid} by considering hybrid inclusions, i.e., differential and difference inclusions to account for continuous and discrete dynamics. The advantage of this framework is that clocks and logical variables can be included into the system description. Hybrid systems with external inputs as in Definition \ref{def:hybrid_system} have explicitly been presented in \cite{sanfelice2016robust}. Note that the value of the state $\boldsymbol{z}_i$ after a jump is denoted by $\hat{\boldsymbol{z}}_i$. This is not a standard convention, but will ease the reading in the upcoming sections.
\begin{definition}\cite{sanfelice2016robust}\label{def:hybrid_system}
A hybrid system is a tuple $\mathcal{H}_i:=(C_i,F_i,D_i,G_i)$ where $C_i$, $D_i$, $F_i$, and $G_i$ are the flow and jump set and the possibly set-valued flow and jump map, respectively. The continuous and discrete dynamics are
\begin{align}\label{eq:hybrid_systems}
\begin{cases}
\dot{\boldsymbol{z}}_i\in F_i(\boldsymbol{z}_i,\boldsymbol{u}_i^{\text{int}},\boldsymbol{u}_i^{\text{ext}}) & \text{for}\;\;\;(\boldsymbol{z}_i,\boldsymbol{u}^{\text{int}}_i,\boldsymbol{u}^{\text{ext}}_i) \in C_i  \\
\hat{\boldsymbol{z}}_i\in G_i(\boldsymbol{z}_i,\boldsymbol{u}^{\text{int}}_i,\boldsymbol{u}^{\text{ext}}_i) & \text{for}\;\;\;(\boldsymbol{z}_i,\boldsymbol{u}^{\text{int}}_i,\boldsymbol{u}^{\text{ext}}_i) \in D_i, 
\end{cases}
\end{align}
where $\boldsymbol{z}_i\in\mathcal{Z}_i$ is a hybrid state with domain $\mathcal{Z}_i$, while $\boldsymbol{u}^{\text{int}}_i\in\mathcal{U}_i^{\text{int}}$ and $\boldsymbol{u}^{\text{ext}}_i \in\mathcal{U}_i^{\text{ext}}$ are internal and external inputs with domains $\mathcal{U}_i^{\text{int}}$ and $\mathcal{U}_i^{\text{ext}}$. Furthermore, let $\mathfrak{H}_i:=\mathcal{Z}_i\times \mathcal{U}_i^{\text{int}}\times \mathcal{U}_i^{\text{ext}}$.
\end{definition}

Solutions to \eqref{eq:hybrid_systems} are parametrized by $(t,j)$, where $t$ indicates continuous flow according to $F_i(\boldsymbol{z}_i,\boldsymbol{u}_i^{\text{int}},\boldsymbol{u}_i^{\text{ext}})$ and $j$ indicates jumps according to $G(\boldsymbol{z}_i,\boldsymbol{u}_i^{\text{int}},\boldsymbol{u}_i^{\text{ext}})$. Hence, a solution is a function $\boldsymbol{z}_i:\mathbb{R}_{\ge 0}\times \mathbb{N} \to \mathcal{Z}_i$ that satisfies \eqref{eq:hybrid_systems} with initial condition $\boldsymbol{z}_i(0,0)$. For a detailed review of the topic, the reader is referred to \cite{goebel2012hybrid}.
\section{Problem Statement}
\label{sec:problem_statement}
In this paper, the following STL fragment is considered:
\begin{subequations}\label{eq:subclass}
\begin{align}
&\hspace{0.12cm}\psi \; ::= \; \top \; | \; \mu \; | \; \neg \mu \; | \; \psi_{(1)} \wedge \psi_{(2)}\label{eq:psi_class}\\
&\hspace{0.15cm}\phi \; ::= \;  G_{[a,b]}\psi \; | \; F_{[a,b]} \psi \label{eq:phi_class}\\
&\theta^{\text{s}_1} ::= \bigwedge_{k=1}^{K} \phi_{(k)} \text{ with } b_{(k)}\le a_{(k+1)} \label{eq:theta1_class} \\
&\theta^{\text{s}_2} ::=  F_{[c_{(1)},d_{(1)}]} \big( \psi_{(1)} \wedge F_{[c_{(2)},d_{(2)}]}(\psi_{(2)} \wedge \hdots ) \big)\label{eq:theta2_class}\\
&\hspace{0.275cm}\theta ::= \theta^{s_1} \; |\; \theta^{s_2},\label{eq:theta_class}
\end{align}
\end{subequations}
where $\mu$ is a predicate and $\psi_{(1)}$, $\psi_{(2)}, \hdots$ are formulas of class $\psi$ given in \eqref{eq:psi_class}, whereas $\phi_{(k)}$ with $k\in\{1,\hdots,K\}$ are formulas of class $\phi$ given in \eqref{eq:phi_class} with corresponding time intervals $[a_{(k)},b_{(k)}]$. Note the use of brackets, e.g. $\psi_{(1)}$, to distinguish from local formulas, e.g., $\psi_{1}$. In this paper, for conjunctions of non-temporal formulas of class $\psi$ given in \eqref{eq:psi_class}, e.g., $\psi:=\psi_{(1)}\wedge\psi_{(2)}$, we approximate the robust semantics, e.g., $\rho^{\psi_{(1)}\wedge\psi_{(2)}}(\boldsymbol{x})$,  by a smooth function.
\begin{assumption}\label{assumption:2}
The robust semantics for a conjunction of $q$ \emph{non-temporal formulas} of class $\psi$ given in \eqref{eq:psi_class}, i.e., $\rho^{\psi_{(1)} \wedge \hdots \wedge \psi_{(q)}}(\boldsymbol{x})$, are approximated by a smooth function as 
\begin{align*}
\rho^{\psi_{(1)} \wedge \hdots \wedge \psi_{(q)}}(\boldsymbol{x})\approx -\ln\Big(\sum_{i=1}^q \exp\big(-\rho^{\psi_{(i)}}(\boldsymbol{x})\big)\Big).
\end{align*}
\end{assumption}
From now on, when writing $\rho^{\psi}(\boldsymbol{x})$, $\rho^{\phi}(\boldsymbol{x},t)$, or $\rho^{\theta}(\boldsymbol{x},t)$ for formulas of class $\psi$, $\phi$, and $\theta$, respectively, we mean the robust semantics including the smooth approximation in Assumption \ref{assumption:2} unless stated otherwise.  This approximation is an under-approximation and preserves the property $(\boldsymbol{x},0)\models\psi$ if $\rho^{\psi}(\boldsymbol{x})>0$ as in \cite{lindemann2017prescribed}. 

The objective in this paper is to consider local formulas of class $\phi$ given in  \eqref{eq:phi_class} that are independently distributed to each agent $v_i\in\mathcal{V}$. The proposed solution can then be extended to local formulas of class $\theta$ given \eqref{eq:theta_class} in the same vein as in \cite{lindemann2017prescribed}. In \cite{lindemann2017prescribed}, a continuous feedback control law for a single agent subject to $\phi$ has been derived, however not considering possible multi-agent couplings as given by $f_i^\text{c}(\boldsymbol{x})$ or formula-agent dependencies. Assume hence that each agent $v_i\in\mathcal{V}$ is subject to a local formula $\phi_i$ of the form \eqref{eq:phi_class}. Two more assumption are needed.

\begin{assumption}\label{assumption:4}
Each formula of class $\psi$ given in \eqref{eq:psi_class} that is contained in \eqref{eq:phi_class} and associated with an agent $v_i$ is: 1) s.t. $\rho^{\psi_i}(\boldsymbol{x}_{\phi_i})$ is concave and 2) well-posed in the sense that $(\boldsymbol{x}_{\phi_i},0)\models \psi_i$ implies $\|\boldsymbol{x}_{\phi_i}(0)\|\le C <\infty$ for some $C\ge 0$.
\end{assumption}
\begin{remark}\label{rem:well-posed}
Part 2) of Assumption  \ref{assumption:4} is not restrictive in practice since $\psi^{\text{Ass.3}}_i:=(\|\boldsymbol{x}_{\phi_i}\|\le C)$, where $C$ is a sufficiently large positive constant, can be combined with the desired $\psi_i$ so that $\psi_i\wedge\psi^{\text{Ass.3}}_i$ is well-posed.
\end{remark}

Next, define the global optimum of $\rsnt{\psi_i}$ as 
\begin{align*}
\rso{\psi_i}:=\sup_{\boldsymbol{x}_{\phi_i}\in\mathbb{R}^{nP(\phi_i)}} \rsnt{\psi_i},
\end{align*} 
which is straightforward to compute due to Assumption \ref{assumption:2} and \ref{assumption:4}. Next, Assumption \ref{assumption:3} guarantees that $\phi_i$ is locally feasible since $\rso{\psi_i}> 0$ implies that $\rsct{\phi_i}{0}> 0$ is possible.
\begin{assumption}\label{assumption:3}
The optimum of $\rsnt{\psi_i}$ is s.t. $\rso{\psi_i}> 0$.
\end{assumption}

The goal is to derive a local control law $\boldsymbol{u}_i(\boldsymbol{x}_{\phi_i},t)$ for each agent $v_i$ such that $r_i\le \rsct{\phi_i}{0}\le \rho^{\text{max}}_i$, where $r_i\in\mathbb{R}$ is a robustness measure, while $\rho^{\text{max}}_i\in\mathbb{R}$ with $r_i<\rho^{\text{max}}_i$ is a robustness delimiter. For this purpose, we look at each dependency cluster separately and distinguish between two cases that are described in the formal problem definition. 
\begin{problem}\label{problem1}
Assume that each agent $v_i$ is subject to a local STL formula $\phi_i$ of the form \eqref{eq:phi_class}, hence inducing the maximal dependency clusters $\bar{\Xi}:=\{\Xi_1,\hdots,\Xi_L\}$ with $L\le M$. For each cluster $\Xi_l$ with $l\in\{1,\hdots,L\}$, derive a control strategy as follows: 
\begin{itemize}
\item Case A) Under the assumption that each agent $v_i,v_j\in\Xi_l$ is subject to the same formula, i.e., $\phi_i=\phi_j$, design a local feedback control law $\boldsymbol{u}_{i}(\boldsymbol{x}_{\phi_i},t)$ such that $0<r_i\le\rsct{\phi_i}{0}\le\rho_i^{\text{max}}$ for all $v_i\in\Xi_l$, which means local satisfaction of $\phi_i$. 
\item Case B) Otherwise, i.e., $\exists v_i,v_j\in\Xi_l$ such that $\phi_i\neq \phi_j$, each agent $v_i\in\Xi_l$ nevertheless initially applies the derived control law $\boldsymbol{u}_{i}(\boldsymbol{x}_{\phi_i},t)$ for Case A. Design a local online detection \& repair scheme for each agent $v_i\in\Xi_l$ such that $r_i\le \rsct{\phi_i}{0}\le \rho_i^{\text{max}}$, where $r_i\in\mathbb{R}$, possibly negative, is maximized up to a precision of $\delta_i>0$ with $\delta_i$ being a design parameter
\end{itemize}
\end{problem}  

If each cluster satisfies the assumption in Case A, i.e., for each $l\in\{1,\hdots,L\}$ it holds that $\phi_i=\phi_j$ for all agents $v_i,v_j\in \Xi_l$,  the proposed solution guarantees global satisfaction of $\{\phi_1,\hdots,\phi_M\}$. If one or more cluster fails to satisfy this assumption, the online detection \& repair scheme in Case B will apply for all agents in these clusters.

\section{Proposed Problem Solution}
\label{sec:inf}

Two types of inter-agent dependencies have been introduced in Section \ref{sec:preliminaries}: dynamic couplings induced by $f_i^\text{c}(\boldsymbol{x})$ and formula-agent dependencies. In the proposed solution to Case A in Problem \ref{problem1}, presented in Section \ref{sec:global_guaran}, it turns out that formula-agent dependencies do not pose any difficulties. Similarly, dynamic couplings only increase the control effort, i.e., $\|\boldsymbol{u}_i(t)\|$.  For Case B, however, both types of dependencies may lead to trajectories that do not locally satisfy the formulas. The proposed solution, introducing an online detection \& repair scheme, is presented in Section \ref{sec:online_repair}.

\subsection{A Prescribed Performance Approach}
We first present the main idea of our work on single-agent systems \cite{lindemann2017prescribed}, which is based on prescribed performance control \cite{bechlioulis2014low} and now extended to multi-agent systems. For a thorough illustration, the reader is referred to \cite{lindemann2017prescribed}. Define the performance function $\gamma_i$ for agent $v_i$ in Definition \ref{def:p} and the transformation function $S$ in Definition \ref{def:SS}. 
\begin{definition}\label{def:p}
The performance function $\gamma_i:\mathbb{R}_{\ge 0}\to\mathbb{R}_{> 0}$ is continuously differentiable, bounded, positive, non-increasing, and given by $\gamma_i(t):=(\gamma_i^0-\gamma_i^\infty)\exp({-l_it})+\gamma_i^\infty$ where $\gamma_i^0, \gamma_i^\infty \in \mathbb{R}_{> 0}$ with $\gamma_i^0\ge \gamma_i^\infty$ and $l_i\in\mathbb{R}_{\ge 0}$.
\end{definition} 
\begin{definition}\label{def:SS}
A transformation function $S:(-1,0)\to\mathbb{R}$ is a strictly increasing function, hence injective and admitting an inverse. In particular, let $S(\xi):=\ln\left(-\frac{\xi+1}{\xi}\right)$.
\end{definition} 

The objective is to synthesize a local feedback control law $\boldsymbol{u}_i(\boldsymbol{x}_{\phi_i},t)$ for formulas $\phi_i$ of the form \eqref{eq:phi_class} such that $r_i \le \rsct{\phi_i}{0} \le \rho^{\text{max}}_i$. Let $\psi_i$ correspond to $\phi_i$ as in $\phi_i:=G_{[a_i,b_i]}\psi_i$ or $\phi_i:=F_{[a_i,b_i]}\psi_i$ and note that $\boldsymbol{x}_{\psi_i}=\boldsymbol{x}_{\phi_i}$ holds by definition. We achieve $r_i\le \rsct{\phi_i}{0}\le \rho^{\text{max}}_i$ by prescribing a temporal behavior to $\rsnttnt{\psi_i}$ through the design parameters $\gamma_i$ and $\rho^{\text{max}}_i$ as 
\begin{align}
 -\gamma_i(t)+\rho^{\text{max}}_i< \rsnttnt{\psi_i}< \rho^{\text{max}}_i. \label{eq:inequality}
\end{align}
Note the use of $\rsnttnt{\psi_i}$ and not $\rsct{\phi_i}{0}$ itself. When $\boldsymbol{x}_{\phi_i}$
is seen as a state, define the one-dimensional error, the normalized error, and the transformed error as \begin{align*}
&\hspace{0.275cm}e_i(\boldsymbol{x}_{\phi_i}):=\rsnt{\psi_i}-\rho^{\text{max}}_i\\
&\xi_i(\boldsymbol{x}_{\phi_i},t):=\frac{e_i(\boldsymbol{x}_{\phi_i})}{\gamma_i(t)}\\
&\hspace{0.005cm}\epsilon_i(\boldsymbol{x}_{\phi_i},t):=S\big(\xi_i(\boldsymbol{x}_{\phi_i},t)\big)=\ln\Big(-\frac{\xi_i(\boldsymbol{x}_{\phi_i},t)+1}{\xi_i(\boldsymbol{x}_{\phi_i},t)} \Big),
\end{align*} 
respectively. Now, \eqref{eq:inequality} can be written as $-\gamma_i(t)< e_i(t)< 0$ where $e_i(t):=e_i(\boldsymbol{x}_{\phi_i}(t))$, which can be further written as $-1< \xi_i(t)< 0$ where $\xi_i(t):=\xi_i(\boldsymbol{x}_{\phi_i}(t),t)$. Applying the transformation function $S$ to $-1< \xi_i(t)< 0$ gives $-\infty<\epsilon_i(t)<\infty$ with $\epsilon_i(t):=\epsilon_i(\boldsymbol{x}_{\phi_i}(t),t)$. If $\epsilon_i(t)$ is bounded for all $t\ge 0$, then inequality \eqref{eq:inequality} holds. This is a consequence of the fact that $S$ admits an inverse. The connection between $\rsnttnt{\psi_i}$ and $\rsct{\phi_i}{0}$ is made by the performance function $\gamma_i$, which needs to be chosen as explained in detail in \cite{lindemann2017prescribed} to obtain $0<r_i\le \rsct{\phi_i}{0}\le \rho^{\text{max}}_i$. If Assumption \ref{assumption:3} holds, then select the parameters
\begin{align}
&\hspace{0.445cm}t^*_i\in\begin{cases} a_i &\text{ if } \phi_i=G_{[a_i,b_i]}\psi_i \\ 
[a_i,b_i] &\text{ if } \phi_i=F_{[a_i,b_i]}\psi_i,
\end{cases}\label{c_t_star}\\
&\hspace{0.1cm}\rho^{\text{max}}_i\in\big(\max\big(0,\rho^{\psi_i}(\boldsymbol{x}_{\phi_i}(0))\big),\rso{\psi_i}\big)\label{rho_max}\\
&\hspace{0.45cm}r_i\in(0,\rho^{\text{max}}_i)\label{r_i}\\
&\hspace{0.385cm}\gamma^{0}_i\in\begin{cases}(\rho^{\text{max}}_i-\rho^{\psi_i}(\boldsymbol{x}_{\phi_i}(0)),\infty) &\text{if } t^*_i>0\\
(\rho^{\text{max}}_i-\rho^{\psi_i}(\boldsymbol{x}_{\phi_i}(0)),\rho^{\text{max}}_i-r_i] &\text{otherwise} \end{cases}\label{eq:g1}\\
&\hspace{0.22cm}\gamma^{\infty}_i\in \Big(0,\min\big(\gamma_i^0,\rho^{\text{max}}_i-r_i\big)\Big]\label{eq:g2}\\
&\hspace{0.5cm}l_i\in \begin{cases}
\mathbb{R}_{\ge 0} & \text{if } -\gamma_i^0+\rho^{\text{max}}_i\ge r_i\\
\frac{-\ln\big(\frac{r_i+\gamma^{\infty}_i-\rho^{\text{max}}_i}{-(\gamma^{0}_i-\gamma^{\infty}_i)}\big)}{t^*_i} & \text{if } -\gamma_i^0+\rho^{\text{max}}_i< r_i\label{eq:g3}
\end{cases}
\end{align}
where it has to hold that $\rho^{\psi_i}(\boldsymbol{x}_{\phi_i}(0))> r_i$ if $t^*_i=0$. The intuition here is that by the choice of $\gamma_i$ it is ensured that $\rho^{\psi_i}(\boldsymbol{x}_{\phi_i}(t))\ge r_i$ for all $t\ge t^*_i$. By the choice of $t^*_i$ it consequently holds that $\rho^{\phi_i}(\boldsymbol{x}_{\phi_i},0)\ge r_i$, i.e., $(\boldsymbol{x}_{\phi_i},0)\models\phi_i$.

\subsection{Global and Local  Satisfaction Guarantees}
\label{sec:global_guaran}
Considering the induced maximal dependency clusters $\bar{\Xi}:=\{\Xi_1,\hdots,\Xi_L\}$, Theorem \ref{theorem:1} provides a global satisfaction guarantee if all clusters satisfy the assumption of Case A in Problem \ref{problem1}, i.e., for each $l\in\{1,\hdots,L\}$ it holds that $\phi_i=\phi_j$ for all $v_i,v_j\in\Xi_l$.

\begin{theorem}\label{theorem:1}
Let each agent $v_i\in\mathcal{V}$ be subject to $\phi_i$ as in \eqref{eq:phi_class}, hence inducing the maximal dependency clusters $\bar{\Xi}:=\{\Xi_1,\hdots,\Xi_L\}$. Assume that for each $\Xi_l\in\bar{\Xi}$ it holds that: for all $v_i,v_j\in \Xi_l$ we have 1) $v_i$ and $v_j$ can communicate, 2) $\phi_i=\phi_j$, and 3) $t^*_i=t^*_j$, $\rho^{\text{max}}_i=\rho^{\text{max}}_j$, $r_i=r_j$, and $\gamma_i=\gamma_j$ are chosen as in \eqref{c_t_star}-\eqref{eq:g3}.   If for each agent $v_i\in\mathcal{V}$ Assumptions \ref{assumption:1}-\ref{assumption:3} hold and each agent $v_i$ applies 
\begin{align}\label{equ:control1}
\boldsymbol{u}_i(\boldsymbol{x}_{\phi_i},t):=-\epsilon_i(\boldsymbol{x}_{\phi_i},t){g_i(\boldsymbol{x}_i)}^T\frac{\partial \rsnt{\psi_i}}{\partial \boldsymbol{x}_i},
\end{align}
then it holds that $0<r_i\le \rsct{\phi_i}{0}\le \rho^{\text{max}}_i$ for all agents $v_i\in \mathcal{V}$, i.e., each agent $v_i$ locally satisfies $\phi_i$, which in turn guarantees global satisfaction of $\{\phi_1,\hdots,\phi_M\}$. All closed-loop signals are well-posed, i.e., continuous and bounded. 

\begin{proof}
In a first step (Step A), we apply Lemma \ref{theorem:sontag1} and show that there exists a maximal solution $\xi_i(t)$ such that $\xi_i(t):=\xi_i(\boldsymbol{x}_{\phi_i}(t),t)\in\Omega_\xi:=(-1,0)$, which is the same as requiring that \eqref{eq:inequality} holds for all $t\in\mathcal{J}:=[0,\tau_{\text{max}})\subseteq \mathbb{R}_{\ge 0}$ and all $v_i\in\Xi_l$ with $l\in\{1,\hdots,L\}$. The second step (Step B) consists of using Lemma \ref{theorem:sontag2} to show that $\tau_{\text{max}}=\infty$, which proves the main result. 

Prior to Step A and B, we state the dynamics of $\epsilon_i$ as
\begin{align}\label{eq:eps_dyn}
\frac{d\epsilon_i}{dt} = \frac{\partial \epsilon_i}{\partial \xi_i} \frac{d\xi_i}{dt}=-\frac{1}{\gamma_i\xi_i(1+\xi_i)}\Big(\frac{\partial \rsnt{\psi_i}}{\partial\boldsymbol{x}}^T\dot{\boldsymbol{x}}-\xi_i\dot{\gamma}_i\Big),
\end{align} 
which can be derived since it holds that $\frac{\partial \epsilon_i}{\partial \xi_i} = -\frac{1}{\xi_i(1+\xi_i)}$ and 
\begin{align}\label{eq:xi_dyn}
\frac{d\xi_i}{dt}=\frac{1}{\gamma_i}\Big(\frac{de_i}{dt}-\xi_i\dot{\gamma}_i\Big).
\end{align} 
Note that by $\frac{d \epsilon_i}{dt}$, $\frac{d \xi_i}{dt}$, and $\frac{d e_i}{dt}$ we here mean the total derivative and that hence $\frac{de_i}{dt} = \frac{\partial e(\boldsymbol{x}_{\phi_i})}{\partial\boldsymbol{x}}^T\dot{\boldsymbol{x}}$ with $\frac{\partial e_i(\boldsymbol{x}_{\phi_i})}{\partial\boldsymbol{x}}=\frac{\partial \rsnt{\psi_i}}{\partial\boldsymbol{x}}$. 

Step A: First, define $\boldsymbol{\xi}:=(
\xi_1 , \xi_2 , \hdots , \xi_M )$ and the stacked vector $\boldsymbol{y}:=(
\boldsymbol{x} , \boldsymbol{\xi})$. Consider the closed-loop system $\dot{\boldsymbol{x}}_i=:H_{\boldsymbol{x}_i}(\boldsymbol{x},\xi_i)$ of agent $v_i$  with  
\begin{align*}
H_{\boldsymbol{x}_i}(\boldsymbol{x},\xi_i)&:=f_i(\boldsymbol{x}_i)+f_i^\text{c}(\boldsymbol{x})\\
&-\ln\Big(-\frac{\xi_i+1}{\xi_i}\Big)g_i(\boldsymbol{x}_i)g_i^T(\boldsymbol{x}_i)\frac{\partial \rsnt{\psi_i}}{\partial\boldsymbol{x}_i}+\boldsymbol{w}_i
\end{align*} 
that is obtained by inserting \eqref{equ:control1} into \eqref{eq:system}. The closed-loop system of all agents is then $\dot{\boldsymbol{x}}=:H_{\boldsymbol{x}}(\boldsymbol{x},\boldsymbol{\xi})$ with
\begin{align*}
H_{\boldsymbol{x}}(\boldsymbol{x},\boldsymbol{\xi}):=\begin{bmatrix} {H_{\boldsymbol{x}_1}(\boldsymbol{x},\xi_1)} & \hdots & {H_{\boldsymbol{x}_M}(\boldsymbol{x},\xi_M)} \end{bmatrix}.
\end{align*} 
Due to \eqref{eq:xi_dyn}, we obtain $\frac{d\xi_i}{dt}=:H_{\xi_i}(\boldsymbol{x},\xi_i,t)$ where
\begin{align*}
H_{\xi_i}(\boldsymbol{x},\xi_i,t):=\frac{1}{\gamma_i(t)}\Big(\frac{\partial \rsnt{\psi_i}}{\partial\boldsymbol{x}} H_{\boldsymbol{x}}(\boldsymbol{x},\boldsymbol{\xi})-\xi_i \dot{\gamma}_i(t)  \Big)
\end{align*} 
expresses the $\xi$-dynamics of agent $v_i$. The $\xi$-dynamics of all agents  are given by $\dot{\boldsymbol{\xi}}=:H_{\boldsymbol{\xi}}(\boldsymbol{x},\boldsymbol{\xi},t)$ with 
\begin{align*}
H_{\boldsymbol{\xi}}(\boldsymbol{x},\boldsymbol{\xi},t):=\begin{bmatrix} H_{\xi_1}(\boldsymbol{x},\xi_1,t) & \hdots & H_{\xi_M}(\boldsymbol{x},\xi_M,t)\end{bmatrix}.
\end{align*}

Using all these definitions, the dynamics of $\boldsymbol{y}$ are finally given by $\dot{\boldsymbol{y}}=:H(\boldsymbol{y},t)$ with
\begin{align*}
H(\boldsymbol{y},t):=\begin{bmatrix}
{H_{\boldsymbol{x}}(\boldsymbol{x},\boldsymbol{\xi})} & {H_{\boldsymbol{\xi}}(\boldsymbol{x},\boldsymbol{\xi},t)}
\end{bmatrix}.
\end{align*} 

Note that $\boldsymbol{x}(0)$ is such that $\xi_i(\boldsymbol{x}_{\phi_i}(0),0)\in\Omega_\xi:=(-1,0)$ holds for all agents $v_i\in\Xi_l$ due to the choice of $\gamma^0_i$. Now define the time-varying and non-empty set 
\begin{align*}
\Omega_{\phi_i}(t)&:=\Big\{\boldsymbol{x}_{\phi_i}\in\mathbb{R}^{nP(\phi_i)}\big|\\
&\hspace{1cm}-1<\xi_i(\boldsymbol{x}_{\phi_i},t)=\frac{\rho^{\psi_i}(\boldsymbol{x}_{\phi_i})-\rho^{\text{max}}_i}{\gamma_i(t)}<0\Big\},
\end{align*} which has the property that for $t_1<t_2$ we have $\Omega_{{\phi_i}}(t_2)\subseteq \Omega_{{\phi_i}}(t_1)$ since $\gamma_i$ is non-increasing in $t$. Note that $\Omega_{\phi_i}(t)$ is bounded due to Assumption \ref{assumption:4} and since $\gamma_i$ is bounded. We remark that $\boldsymbol{x}_{\phi_i}(0)\in\Omega_{{\phi_i}}(0)$. Due to \cite[Proposition 1.4.4]{aubin2009set}, the following holds: if a function is continuous, then the inverse image of an open set under this function is open. By defining $\xi_i^0(\boldsymbol{x}_{\phi_i}):=\xi_i(\boldsymbol{x}_{\phi_i},0)$, it holds that inv$\big(\xi_i^0(\Omega_\xi)\big)=\Omega_{{\phi_i}}(0)$ is open. Note therefore that $\rho^{\psi_i}(\boldsymbol{x}_{\phi_i})$ is a continuously differentiable function due to Assumption \ref{assumption:2}. Next, select $v_{i_l}\in\Xi_l$ for each $l\in\{1,\hdots,L\}$ and define 
\begin{align*}
\Omega_{\boldsymbol{x}}&:=\Omega_{{\phi_{i_1}}}(0)\times \hdots \times \Omega_{{\phi_{i_L}}}(0)\subset\mathbb{R}^{nM},\\
\Omega_{\boldsymbol{\xi}}&:={\Omega_\xi}\times \hdots \times {\Omega_\xi}\subset \mathbb{R}^M,
\end{align*}
and the open, non-empty, and bounded set 
\begin{align*}
\Omega_{\boldsymbol{y}}:=\Omega_{\boldsymbol{x}} \times \Omega_{\boldsymbol{\xi}}\subset \mathbb{R}^{(n+1)M}
\end{align*} 
where it holds that $\boldsymbol{y}(0)=\begin{bmatrix}
{\boldsymbol{x}(0)} & {\boldsymbol{\xi}(0)}
\end{bmatrix}\in\Omega_{\boldsymbol{y}}$. 

Next, we check the conditions in Lemma \ref{theorem:sontag1} for the initial value problem $\dot{\boldsymbol{y}}=H(\boldsymbol{y},t)$ with $\boldsymbol{y}(0)\in\Omega_{\boldsymbol{y}}$ and $H(\boldsymbol{y},t):\Omega_{\boldsymbol{y}}\times \mathbb{R}_{\ge 0} \to \mathbb{R}^{(n+1)M}$: 1) $H(\boldsymbol{y},t)$ is locally Lipschitz continuous on $\boldsymbol{y}$ since $f_i(\boldsymbol{x}_i)$, $f_i^\text{c}(\boldsymbol{x})$, $g_i(\boldsymbol{x}_i)$, and $\epsilon_i=\ln\big(-\frac{\xi_i+1}{\xi_i}\big)$ are locally Lipschitz continuous on $\boldsymbol{y}$ for each $t\in \mathbb{R}_{\ge 0}$. This also holds for $\frac{\partial \rsnt{\psi_i}}{\partial\boldsymbol{x}_i}$ due to Assumption \ref{assumption:2}. 2) $H(\boldsymbol{y},t)$ is continuous on $t$ for each fixed $\boldsymbol{y}\in\Omega_{\boldsymbol{y}}$ due to continuity of $\gamma_i$ and $\dot{\gamma}_i$. As a result of Lemma \ref{theorem:sontag1}, there exists a maximal solution with $\boldsymbol{y}(t)\in\Omega_{\boldsymbol{y}}$ for all $t\in\mathcal{J}:=[0,\tau_{\text{max}})\subseteq \mathbb{R}_{\ge 0}$ and $\tau_{\text{max}}>0$. Consequently, there exist $\boldsymbol{\xi}(t)\in\Omega_{\boldsymbol{\xi}}$ and $\boldsymbol{x}(t)\in\Omega_{\boldsymbol{x}}$ for all $t \in \mathcal{J}$. 

Step B: From Step A) we have $\boldsymbol{y}(t)\in\Omega_{\boldsymbol{y}}$ for all $t\in\mathcal{J}:=[0,\tau_{\text{max}})$. Now, we show that $\tau_{\text{max}}=\infty$ by contradiction of Lemma \ref{theorem:sontag2}. Therefore, assume $\tau_{\text{max}}<\infty$.  

The key observation to be made is that $\xi_i(\boldsymbol{x}_{\phi_i},t)=\xi_j(\boldsymbol{x}_{\phi_j},t)$, $\epsilon_i(\boldsymbol{x}_{\phi_i},t)=\epsilon_j(\boldsymbol{x}_{\phi_j},t)$, and $\rsnt{\psi_i}=\rsnt{\psi_j}$ for all agents $v_i,v_j\in\Xi_l$. This follows since $\boldsymbol{x}_{\phi_i}=\boldsymbol{x}_{\phi_j}$ (recall that $\phi_i=\phi_j$) and since $\rho^{\text{max}}_i=\rho^{\text{max}}_j$ and $\gamma_i=\gamma_j$ holds by assumption. We now show that $\epsilon_i(t)$ is bounded for all $t\in\mathbb{R}_{\ge 0}$ and then it consequently follows that $\epsilon_j(t)$ is bounded for all other agents $v_j\in\Xi_l\setminus\{v_i\}$. Since the clusters are maximal, i.e., no formula-agent dependencies between clusters exist, we can deduce the same result for the other clusters. Consider the Lyapunov function candidate $V(\epsilon_i):=\frac{1}{2}\epsilon_i\epsilon_i$ and define $\dot{V}(\epsilon_i):=\frac{\partial V}{\partial \epsilon_i}\frac{d \epsilon_i}{dt}$. We will now show that $\dot{V}(\epsilon_i)\le 0$ if $|\epsilon_i|$ is bigger than some positive constant, which ensures that $\epsilon_i(t)$ will remain in a compact set. By using \eqref{eq:eps_dyn}, it follows
\begin{align*}
\dot{V}(\epsilon_i)&=\epsilon_i\frac{d\epsilon_i}{dt}=\epsilon_i\Big(-\frac{1}{\gamma_i\xi_i(1+\xi_i)}\big(\frac{\partial \rsnt{\psi_i}}{\partial \boldsymbol{x}}^T\dot{\boldsymbol{x}}-\xi_i\dot{\gamma}_i\big)\Big)
\end{align*} 
Define $\alpha_i(t):=-\frac{1}{\gamma_i\xi_i(1+\xi_i)}$ which satisfies $\alpha_i(t)\in[\frac{4}{\gamma_i^0},\infty)$ for all $t\in \mathcal{J}$. This follows since $\frac{4}{\gamma^0_i}\le -\frac{1}{\gamma^0_i\xi_i(1+\xi_i)}\le-\frac{1}{\gamma_i\xi_i(1+\xi_i)}\le-\frac{1}{\gamma_i^\infty\xi_i(1+\xi_i)}< \infty$ for $\xi_i\in\Omega_\xi$. It can further be derived that
\begin{align}\label{eq:V_raw111}
\dot{V}(\epsilon_i)\le \epsilon_i\alpha_i\frac{\partial \rsnt{\psi_i}}{\partial \boldsymbol{x}}^T  H_{\boldsymbol{x}}(\boldsymbol{x},\boldsymbol{\xi}) +|\epsilon_i|\alpha_i k_i 
\end{align}
where $\dot{\boldsymbol{x}}:=H_{\boldsymbol{x}}(\boldsymbol{x},\boldsymbol{\xi})$ as defined previously and $0\le |\xi_i\dot{\gamma}_i|\le k_i<\infty$ for a positive constant $k_i$. This follows since $\xi_i(t)\in\Omega_\xi$ for all $t\in\mathcal{J}$ and $\dot{\gamma}_i$ is bounded by definition. The term $\frac{\partial \rsnt{\psi_i}}{\partial \boldsymbol{x}}^T  H_{\boldsymbol{x}}(\boldsymbol{x},\boldsymbol{\xi})$ represents the couplings among the agents and can be written as 
\begin{align}\label{eq:sum_H}
\frac{\partial \rsnt{\psi_i}}{\partial \boldsymbol{x}}^T  H_{\boldsymbol{x}}(\boldsymbol{x},\boldsymbol{\xi})=
\sum_{v_j\in\Xi_l}\frac{\partial \rsntj{\psi_j}}{\partial \boldsymbol{x}_j}^T  H_{\boldsymbol{x}_j}(\boldsymbol{x},\xi_j).
\end{align} 
Plugging \eqref{eq:sum_H} into \eqref{eq:V_raw111} results in 
\begin{align}\label{eq:V_raw11}
\dot{V}(\epsilon_i)\le \epsilon_i\alpha_i\sum_{v_j\in\Xi_l}\frac{\partial \rsntj{\psi_j}}{\partial \boldsymbol{x}_j}^T  H_{\boldsymbol{x}_j}(\boldsymbol{x},\xi_j) +|\epsilon_i|\alpha_i k_i.
\end{align}
Inserting \eqref{eq:system} and \eqref{equ:control1} into $\epsilon_i\frac{\partial \rsntj{\psi_j}}{\partial \boldsymbol{x}_j}^T  H_{\boldsymbol{x}_j}(\boldsymbol{x},\xi_j)$ first, this term can in a second step be upper bounded as follows 
\begin{align*}
&\epsilon_i\frac{\partial \rsntj{\psi_j}}{\partial \boldsymbol{x}_j}^T  H_{\boldsymbol{x}_j}(\boldsymbol{x},\xi_j)=\epsilon_i\frac{\partial \rsntj{\psi_j}}{\partial \boldsymbol{x}_j}^T \Big(f_j(\boldsymbol{x}_j)+f_j^\text{c}(\boldsymbol{x})\\
&\hspace{3.5cm}-\epsilon_jg_j(\boldsymbol{x}_j)g_j^T(\boldsymbol{x}_j)\frac{\partial \rsntj{\psi_j}}{\partial\boldsymbol{x}_j}+\boldsymbol{w}_j\Big)\\
&\hspace{3.75cm}\le |\epsilon_i|M_j-|\epsilon_i|^2\lambda_j J_j
\end{align*}
where $\epsilon_i=\epsilon_j$ as remarked previously since $v_i,v_j\in\Xi_l$. Furthermore, $\lambda_j>0$ is the positive minimum eigenvalue of $g_j(\boldsymbol{x}_j)g_j^T(\boldsymbol{x}_j)$ according to Assumption \ref{assumption:1}, and $\|\frac{\partial \rsntj{\psi_j}}{\partial \boldsymbol{x}_j}^T \big(f_j(\boldsymbol{x}_j)+f_j^\text{c}(\boldsymbol{x})+\boldsymbol{w}_j\big)\|\le M_j<\infty$ due to continuity of $\frac{\partial \rsntj{\psi_j}}{\partial \boldsymbol{x}_j}$, $f_j(\boldsymbol{x}_j)$, and $f_j^\text{c}(\boldsymbol{x})$, the extreme value theorem and the fact that $\Omega_{\boldsymbol{x}}$ and $\mathcal{W}_i$ are bounded.  Note therefore that the extreme value theorem guarantees that a continuous function on a compact set is bounded and that the above functions are continuous on $\text{cl}(\Omega_{\boldsymbol{x}})$, where $\text{cl}$ denotes the closure of a set. The lower bound $J_j\in\mathbb{R}_{\ge 0}$ arises naturally due to the norm operator as $0\le J_j\le  \|\frac{\partial \rsntj{\psi_j}}{\partial\boldsymbol{x}_j})\|^2<\infty$.  Equation \eqref{eq:V_raw11} can now be upper bounded as follows
\begin{align}
\dot{V}(\epsilon_i)\le \alpha_i|\epsilon_i|\big(\hat{M}_i -|\epsilon_i|\hat{J}_i\big)\
\end{align}
where $\hat{M}_i:=\sum_{v_j\in\Xi_l}M_j+ k_i$ and $\hat{J}_i:=\sum_{v_j\in\Xi_l}\lambda_jJ_j$. Note that $\hat{J}_i>0$ since $\frac{\partial \rho^{\psi_i}(\boldsymbol{x}_{\phi_i})}{\partial\boldsymbol{x}_{\phi_i}}=0$ if and only if $\rho^{\psi_i}(\boldsymbol{x}_{\phi_i})=\rso{\psi_i}$, which is excluded  since \eqref{eq:inequality} holds for all $t\in\mathcal{J}$ and we selected $\rho^{\text{max}}_i<\rso{\psi_i}$. Recall that $\rsntj{\psi_j}$ in $\|\frac{\partial \rsntj{\psi_j}}{\partial\boldsymbol{x}_j}\|^2$ is concave due to Assumptions \ref{assumption:2} and \ref{assumption:4}. In other words, at least one $J_j$ in $v_j\in\Xi_l$ is greater than zero.

It holds that $\dot{V}(\epsilon_i)\le 0$ if  $\frac{\hat{M}_i}{\hat{J}_i}\le  |\epsilon_i|$. We can conclude that $|\epsilon_i|$ will be upper bounded due to the level sets of $V$ as 
\begin{align*}
|\epsilon_i(t)|\le \max\Big(|\epsilon_i(0)|,\frac{\hat{M}_i}{\hat{J}_i}\Big),
\end{align*} 
which leads to the conclusion that $\epsilon_i(t)$ is upper and lower bounded by some constants $\epsilon_i^u$ and $\epsilon_i^l$, respectively. In other words, it holds that $\epsilon_i^l\le \epsilon_i(t)\le \epsilon_i^u$ for all $t\in\mathcal{J}$. By using the inverse of $S$ and defining $\xi_i^l:=-\frac{1}{\exp(\epsilon_i^l+1)}$ and $\xi_i^u:=-\frac{1}{\exp(\epsilon_i^u+1)}$, $\xi_i(t)$ is bounded by 
$-1<\xi_i^l\le \xi_i(t) \le \xi_i^u<0$, which translates to 
\begin{align*}
\xi_i(t)\in\Omega_{\xi_i}^\prime:= [\xi_i^l,\xi_i^u]\subset \Omega_\xi
\end{align*} 
for all $t\in \mathcal{J}$. Recall that $\xi_i\big(\boldsymbol{x}_{\phi_i},t\big)=\frac{\rsnt{\psi_i}-\rho^{\text{max}}_i}{\gamma_i(t)}$ and note the following: if $\xi_i(t)$ evolves in a compact set, then $\rsnttm{\psi_i}$ will evolve in a compact set $\Omega_{\rho^{\psi_i}}^\prime:=[\rho_i^l,\rho_i^u]$ for some constants $\rho_i^l$ and $\rho_i^u$. Again, due to \cite[Proposition 1.4.4]{aubin2009set} it holds that the inverse image 
\begin{small}
\begin{align*}
\Omega_{{\phi_i}}^\prime:=\text{inv}\big(\rho^{\psi_i}(\Omega_{\rho^{\psi_i}}^\prime)\big)=\{\boldsymbol{x}_{\phi_i}\in \Omega_{{\phi_i}}(0)|\rho_i^l\le \rsnt{\psi_i}\le \rho_i^u\}
\end{align*} 
\end{small}is closed and also bounded since it is a subset of $\Omega_{{\phi_i}}$. Select $v_{i_l}\in\Xi_l$ for each $l\in\{1,\hdots,L\}$. It can be concluded that $\boldsymbol{x}_{\phi_{i_l}}(t)$ evolves in a compact set, i.e., $\boldsymbol{x}_{\phi_{i_l}}(t)\in \Omega_{{\phi_{i_l}}}^\prime\subset\Omega_{{\phi_{i_l}}}(0)$ for all $t\in \mathcal{J}$ and all $v_{i_l}$. Next, define 
\begin{align*}
\Omega_{\boldsymbol{x}}^\prime&:= \Omega_{{\phi_{i_1}}}^\prime\times \hdots \times \Omega_{{\phi_{i_L}}}^\prime\subset\mathbb{R}^{nM}\\
\Omega_{\boldsymbol{\xi}}&:=\Omega_{\xi_1}^\prime\times \hdots \times \Omega_{\xi_M}^\prime\subset\mathbb{R}^M,
\end{align*} 
and the compact set 
\begin{align*}
\Omega_{\boldsymbol{y}}^\prime:=\Omega_{\boldsymbol{x}}^\prime \times \Omega_{\boldsymbol{\xi}}^\prime\subset\mathbb{R}^{(n+1)M}
\end{align*} 
for which it holds that $\boldsymbol{y}(t)\in\Omega_{\boldsymbol{y}}^\prime$ for all $t\in \mathcal{J}$. It is also true that $\Omega_{\boldsymbol{y}}^\prime\subset \Omega_{\boldsymbol{y}}$ by which it follows that there is no $t\in \mathcal{J}:=[0,\tau_{\text{max}})$ such that $\boldsymbol{y}(t)\notin \Omega_{\boldsymbol{y}}^\prime$. By contradiction of Lemma \ref{theorem:sontag2} it holds that $\tau_{\text{max}}=\infty$, i.e., $\mathcal{J}=\mathbb{R}_{\ge 0}$. This in turn says that \eqref{eq:inequality} holds for all agents $v_i\in\mathcal{V}$ and for all $t\in \mathbb{R}_{\ge 0}$. By the choice of $\rho^{\text{max}}_i$, $r_i$, and $\gamma_i$ as in \eqref{rho_max}-\eqref{eq:g3} and \cite[Theorem~2]{lindemann2017prescribed}, it then holds that $\phi_i$ is locally satisfied for each agent $v_i\in\mathcal{V}$. 

The control law $\boldsymbol{u}_i(\boldsymbol{x}_i,t)$ is well-posed, i.e., continuous and bounded, because $\rsntnc{\psi_i}$ is approximated by a smooth function, while $\epsilon_i(\boldsymbol{x}_i,t)$ and $g_i(\boldsymbol{x}_i)$ are continuous. Furthermore, $\gamma_i$ is continuous with $0<\gamma(t)<\infty$. Due to the extreme value theorem, these functions are also bounded. It follows that all closed-loop signals are well-posed.
\end{proof}
\end{theorem}

If $L=M$, i.e., each agent $v_i\in\mathcal{V}$ is subject to a non-collaborative formula $\phi_i$, Theorem \ref{theorem:1} trivially applies since no formula dependencies among agents exist. Recall that dynamic couplings induced by $f^\text{c}_i(\boldsymbol{x})$ may still be present.

%

For the next result, a stronger assumption on the dynamic couplings $f_i^\text{c}(\boldsymbol{x})$ is needed.
\begin{assumption}\label{ass:5}
	The function $f_i^\text{c}:\mathbb{R}^{nM}\to\mathbb{R}^n$ is bounded.
\end{assumption}

Now consider a formula $\phi$ of the form \eqref{eq:phi_class} and assume that each $v_i\in\mathcal{V}_{\phi}$ is subject to $\phi_i:=\phi$. Then, Theorem \ref{theorem:2} guarantees satisfaction of $\phi$ if all agents $v_i\in\mathcal{V}_{\phi}$ collaborate.

\begin{theorem}\label{theorem:2}
Let each agent $v_i\in\mathcal{V}$ satisfy Assumption \ref{assumption:1} and \ref{ass:5}. Consider a formula $\phi$ as in \eqref{eq:phi_class}  and let each agent $v_i\in\mathcal{V}_{\phi}$ be subject to $\phi_i:=\phi$. Assume that for all $v_i,v_j\in \mathcal{V}_{\phi}$ it holds that: 1) $v_i$ and $v_j$ can communicate and 2) $t^*_i=t^*_j$, $\rho^{\text{max}}_i=\rho^{\text{max}}_j$, $r_i=r_j$, and $\gamma_i=\gamma_j$ are chosen as in \eqref{c_t_star}-\eqref{eq:g3}. Assume further that all agents $v_k\in\mathcal{V}\setminus\mathcal{V}_{\phi}$ apply a control law $\boldsymbol{u}^\prime_k$ such that $\boldsymbol{x}_k$ remains in a compact set $\Omega_k^\prime$.  If for each agent $v_i\in\mathcal{V}_{\phi}$ Assumptions \ref{assumption:2}-\ref{assumption:3} hold and each $v_i\in\mathcal{V}_{\phi}$ applies \eqref{equ:control1}, then it holds that $0<r:=r_i\le \rsxt{\phi}{0}\le \rho^{\text{max}}_i=:\rho^{\text{max}}$, i.e., $(\boldsymbol{x}_\phi,0)\models \phi$. All closed-loop signals are well-posed.
	
	\begin{proof}
	The proof is similar to the proof in Theorem \ref{theorem:1} and is provided in the appendix.
					\end{proof}
	\end{theorem}
	
The assumption of $\boldsymbol{u}^\prime_k$ is not restrictive and excludes finite escape time. For instance, if Assumption \ref{ass:5} holds and $\dot{\boldsymbol{x}}_k:=f_k(\boldsymbol{x}_k)$ is asymptotically stable, then the feedback control law $\boldsymbol{u}_k^\prime(\boldsymbol{x}_k):=-{g_k(\boldsymbol{x}_k)}^T\boldsymbol{x}_k$ keeps the state $\boldsymbol{x}_k$ in a compact set. If all agents $v_i\in\mathcal{V}_{\phi}$ apply the control law \eqref{equ:control1} under the conditions in Theorem \ref{theorem:2} to satisfy $\phi$, we refer to this as \emph{collaborative control} in the remainder.  Theorem \ref{theorem:2} has further implications with respect to Case A in Problem \ref{problem1}. Consider again the induced maximal dependency clusters $\bar{\Xi}:=\{\Xi_1,\hdots,\Xi_L\}$. Assume that the cluster $\Xi_l$ with $l\in\{1,\hdots,L\}$ satisfies the assumption of Case A, while there exists another cluster $\Xi_m$ with $m\neq l$ such that $\Xi_m$ does not satisfy this assumption. In other words, for all $v_i,v_j\in\Xi_l$ it holds that $\phi_i=\phi_j$, while $\exists v_i,v_j\in\Xi_m$ with $m\neq l$ such that $\phi_i\neq \phi_j$. Consequently, Theorem \ref{theorem:2} guarantees local satisfaction of $\phi_i$ for all $v_i\in\Xi_l$ without considering task satisfaction of agents in $\mathcal{V}\setminus\Xi_l$.

Note that Assumption \ref{assumption:3} in Theorem \ref{theorem:2} restricts the formula $\phi$ to be locally feasible. However, this assumption can  be relaxed at the expense of not locally satisfying $\phi_i:=\phi$ and instead finding a, possibly least violating, solution by relaxing $r_i$ and $\rho^{\text{max}}_i$. Recall that $\rsct{\phi_i}{0}\ge r_i$ with $r_i< 0$ does not imply local satisfaction of $\phi_i$. 
\begin{corollary}\label{corollary:2}
Assume that all assumptions of Theorem \ref{theorem:1} hold for each agent $v_i\in \mathcal{V}$ except for Assumption \ref{assumption:3} and the choice of $\rho^{\text{max}}_i$ and $r_i$. If  instead $\rho^{\text{max}}_i\in (\rsntont{\psi_i},\rso{\psi_i})$ and $r_i\in (-\infty,\rho^{\text{max}}_i)$, then it holds that $r_i\le \rsct{\phi_i}{0}\le\rho^{\text{max}}_i$ for all agents $v_i\in \mathcal{V}$. 
	
	\begin{proof}
		Follows the same line of proof as in Theorem \ref{theorem:1} and \ref{theorem:2}. Note that it has already been stated in \cite{lindemann2017prescribed} that  $r_i\le\rsct{\phi_i}{0}\le\rho^{\text{max}}_i$ follows from \eqref{eq:inequality} by the choice of $\gamma_i$. Therefore, it consequently holds that $r_i$ can be chosen negative as long as $r_i<\rho^{\text{max}}_i<\rso{\psi_i}$.
	\end{proof}
\end{corollary}

\subsection{An Online Detection \& Repair Scheme}
\label{sec:online_repair}

Assume now that the cluster $\Xi_l$ with $l\in\{1,\hdots,L\}$ may not satisfy the assumption of Case A in Problem \ref{problem1}. We propose that each agent $v_i\in\Xi_l$ initially applies the control law \eqref{equ:control1} with parameters as in \eqref{c_t_star}-\eqref{eq:g3}.   The control law \eqref{equ:control1} consists of two components, one determining the control strength and one the control direction; $\epsilon_i(\boldsymbol{x}_{\phi_i},t)$ determines the control strength. The closer $\xi_i(\boldsymbol{x}_{\phi_i},t)$ gets to $\Omega_\xi:=\{-1,0\}$, i.e., the funnel boundary, the bigger gets $\epsilon_i(\boldsymbol{x}_{\phi_i},t)$ and consequently also $\|\boldsymbol{u}(\boldsymbol{x}_{\phi_i},t)\|$. Note that $\|\boldsymbol{u}(\boldsymbol{x}_{\phi_i},t)\|\to\infty$  as $\xi_i(\boldsymbol{x}_{\phi_i},t)\to \Omega_\xi$. The control direction is determined by $\frac{\partial \rsnt{\psi_i}}{\partial \boldsymbol{x}_i}$, i.e., in which direction control action should mainly happen. In summary, the control law always steers in the direction away from the funnel boundary, and the control effort increases close to the funnel boundary. We reason that applying the control law \eqref{equ:control1} is hence a good initial choice such that $ \phi_i$ will be locally satisfied if the participating agents $\mathcal{V}_{\phi_i}\setminus \{v_i\}$ behave reasonably. The resulting trajectory $\boldsymbol{x}_{\phi_i}$ may, however, hit the funnel boundary, i.e., $\xi_i(\boldsymbol{x}_{\phi_i},t)=\{-1,0\}$, and lead to critical events. 

%

\begin{example}\label{ex:1}
Consider three agents $v_1$, $v_2$, and $v_3$. Agent $v_2$ is subject to the formula $\phi_2:=F_{[5,15]}(\|\boldsymbol{x}_2-\begin{bmatrix}
90 & 90
\end{bmatrix}\|\le 5)$, while agent $v_3$ is subject to $\phi_3:=F_{[5,15]}(\|\boldsymbol{x}_3-\begin{bmatrix}
90 & 10
\end{bmatrix}\|\le 5)$, i.e., both agents are subject to non-collaborative formulas. Agent $v_1$ is subject to the collaborative formula $\phi_1:=G_{[0,15]}(\|\boldsymbol{x}_1-\boldsymbol{x}_2\|\le 10 \wedge \|\boldsymbol{x}_1-\boldsymbol{x}_3\|\le 10)$. Note that the set of formulas $\{\phi_1,\phi_2,\phi_3\}$ is not globally feasible, although each formula itself is locally feasible. Under \eqref{equ:control1}, agents $v_2$ and $v_3$ move to $\begin{bmatrix}
90 & 90
\end{bmatrix}$ and $\begin{bmatrix}
90 & 10
\end{bmatrix}$, respectively. Agent $v_1$ can consequently not satisfy $\phi_1$ and only decrease the robustness such that $r_i<0$ to achieve $r_i\le \rsct{\phi_i}{0}\le 0$ similar to Corollary \ref{corollary:2}.
\end{example}

In Example \ref{ex:1}, the set of local formulas is globally infeasible. However, even if the set $\{\phi_1,\hdots,\phi_M\}$ is globally feasible, there are reasons why the resulting trajectory may not globally satisfy $\{\phi_1,\hdots,\phi_M\}$  as illustrated next.
\begin{example}\label{ex:2}
Consider two agents $v_4$ and $v_5$ with $\phi_4:=F_{[5,10]}(\|\boldsymbol{x}_4-\boldsymbol{x}_5\|\le 10 \wedge \|\boldsymbol{x}_4-\begin{bmatrix} 50 & 70 \end{bmatrix}\|\le 10)$ (collaborative formula) and  $\phi_5:=F_{[5,15]}(\|\boldsymbol{x}_5-\begin{bmatrix}
10 & 10
\end{bmatrix}\|\le 5)$ (non-collaborative formula). Under \eqref{equ:control1}, agent $v_5$ moves to $\begin{bmatrix}
10 & 10
\end{bmatrix}$ by at latest 15 time units. However, agent $v_4$ is forced to move to $\begin{bmatrix} 50 & 70 \end{bmatrix}$ and be close to agent $v_5$ by at latest 10 time units. This may lead to critical events where \eqref{eq:inequality} is violated for agent $v_4$. If agent $v_5$ cooperates, it can first help to locally satisfy $\phi_4$, e.g., by using \emph{collaborative control} as in Theorem \ref{theorem:2}, and locally satisfy $\phi_5$ afterwards.
\end{example}

To overcome these potential problems, we propose an online detection \& repair scheme by using a local hybrid system $\mathcal{H}_i:=(C_i,F_i,D_i,G_i)$ for each agent $v_i\in\Xi_l$. We \emph{detect} critical events that may lead to trajectories that do not locally satisfy $\phi_i$. Then, agent $v_i$ tries to locally \emph{repair} the funnel, i.e., the design parameters $t^*_i$, $\rho^{\text{max}}_i$, $r_i$, and $\gamma_i$, in a first stage. If this is not successful, \emph{collaborative control} as in Theorem \ref{theorem:2} will be considered in a second stage (Example \ref{ex:2}). If \emph{collaborative control} is not applicable, $r_i$ is successively decreased by $\delta_i>0$ in the third stage (Example \ref{ex:1}), where $\delta_i$ is a design parameter. The jump set $D_i$ will \emph{detect} critical events, while the jump map $G_i$ will take \emph{repair} actions. 

Let $\boldsymbol{p}_i^\gamma:=\begin{bmatrix} \gamma_i^0 & \gamma_i^\infty & l_i\end{bmatrix}$ and  $\boldsymbol{p}^\text{f}_i:=\begin{bmatrix} t^*_i  & \rho^{\text{max}}_i & r_i & \boldsymbol{p}_i^\gamma \end{bmatrix}$ contain the parameters that define \eqref{eq:inequality}, and let $\boldsymbol{p}_i^\text{r}:=\begin{bmatrix}\mathfrak{n}_i & \mathfrak{c}_i  \end{bmatrix}$; $\mathfrak{n}_i$ indicates the number of repair attempts in the first repair stage, while $\mathfrak{c}_i$ is used in the second repair stage ($\mathfrak{c}$ for $\mathfrak{c}$ollaborative). If $\mathfrak{c}_i\in\{1,\hdots,M\}$, \emph{collaborative control} as in Theorem \ref{theorem:2} is used to collaboratively satisfy $\phi_{\mathfrak{c}_i}$. If $\mathfrak{c}_i=0$, then agent $v_i$ tries to locally satisfy $\phi_i$ by itself and if $\mathfrak{c}_i=-1$, then agent $v_i$ is free, i.e., not subject to a task. We define the hybrid state as $\boldsymbol{z}_i:=\begin{bmatrix}  \boldsymbol{x}_i & t_i  & {\boldsymbol{p}^\text{f}_i} & {\boldsymbol{p}_i^\text{r}} \end{bmatrix}\in \mathcal{Z}_i$,
where $t_i$ is a clock, $\mathcal{Z}_i:=  \mathbb{R}^n \times \mathbb{R}_{\ge 0} \times \mathbb{R}_{\ge 0}^{6} \times \mathbb{Z}^2$ and $\boldsymbol{z}_i(0,0):=\begin{bmatrix}  {\boldsymbol{x}_i(0)} & 0 & {\boldsymbol{p}^{\text{f}}_i}(0) & \boldsymbol{0}_2 \end{bmatrix}$ with $\mathbb{Z}$ being the set of integers. The elements in $\boldsymbol{p}^{\text{f}}_i(0)$ are as chosen according to \eqref{c_t_star}-\eqref{eq:g3}. Additionally, we choose $\boldsymbol{p}^{\text{f}}_i(0)=\boldsymbol{p}^{\text{f}}_j(0)$ if Case A holds for all agents $v_i,v_j\in \Xi_l$. Next, define 
\begin{align*}
\boldsymbol{u}_i^\text{int}=\begin{cases}\boldsymbol{0}_{m_i} &\text{if } \mathfrak{c}_i=-1\\-\epsilon_i(\boldsymbol{x}_{\phi_i},t_i){g_i(\boldsymbol{x}_i)}^T\frac{\partial \rsnt{\psi_i}}{\partial \boldsymbol{x}_i} & \text{if } \mathfrak{c}_i=0 \\
-\epsilon_{\mathfrak{c}_i}(\boldsymbol{x}_{\phi_{\mathfrak{c}_i}},t_i){g_i(\boldsymbol{x}_i)}^T\frac{\partial \rsntk{\psi_{\mathfrak{c}_i}}}{\partial \boldsymbol{x}_i} &\text{if } \mathfrak{c}_i>0\end{cases}
\end{align*} 
so that the flow map can be written as
\begin{align*}
F_i(\boldsymbol{x}_i,\boldsymbol{u}^{\text{int}}_i,\boldsymbol{u}^{\text{ext}}_i)&:=\\
&\hspace{-1cm}\begin{bmatrix}
 f_i(\boldsymbol{x}_i)+f_i^\text{c}(\boldsymbol{x})+g_i(\boldsymbol{x}_{i})\boldsymbol{u}_{i}^\text{int}+\boldsymbol{w}_i & 1 & {\boldsymbol{0}_6} & {\boldsymbol{0}_2}
\end{bmatrix}.
\end{align*}

External inputs are $\boldsymbol{w}_i$ and $\boldsymbol{x}_{i}^{\text{ext}}$. By assuming $v_i\in\Xi_{l}$, we define $\boldsymbol{\mathfrak{c}}_i^{\text{ext}}:=\begin{bmatrix}
\mathfrak{c}_{j_1} & \hdots & \mathfrak{c}_{j_{|\Xi_l|-1}}\end{bmatrix}$ and $\boldsymbol{p}^\text{f,\text{ext}}_i:=\begin{bmatrix}
\boldsymbol{p}^\text{f}_{j_1} & \hdots & \boldsymbol{p}^\text{f}_{j_{|\Xi_l|-1}}
\end{bmatrix}$ such that $v_{j_1},\hdots,v_{j_{|\Xi_l|-1}}\in\Xi_l\setminus \{v_i\}$. Note that $\boldsymbol{\mathfrak{c}}_i^{\text{ext}}$ and $\boldsymbol{p}^\text{f,\text{ext}}_i$ contain states of all agents in the same dependency cluster $\Xi_l$. Ultimately, define the external input as $\boldsymbol{u}^{\text{ext}}_i:=\begin{bmatrix}
{\boldsymbol{w}_i} & {\boldsymbol{x}_{i}^{\text{ext}}}  &  \boldsymbol{\mathfrak{c}}_i^{\text{ext}} & \boldsymbol{p}^\text{f,\text{ext}}_i
\end{bmatrix}$.

The set $\mathcal{D}_i^\prime$ is used to \emph{detect} a critical event when the funnel in \eqref{eq:inequality} is violated, i.e., when $\xi_i(t_i)\notin \Omega_\xi:=(-1,0)$.
\begin{align*}
\mathcal{D}_i^\prime:= \{(\boldsymbol{z}_i,\boldsymbol{u}^{\text{int}}_i,\boldsymbol{u}^{\text{ext}}_i)\in\mathfrak{H}_i|\xi_i(t_i)\in\{-1,0\},\; \mathfrak{c}_i=0 \}.
\end{align*}

\begin{remark}
Note that $\xi_i(t_i)\in\{-1,0\}$ implies $\epsilon_i(t_i)\to \infty$ and therefore $\boldsymbol{u}_i(t_i)\to\infty$. In practice, the input will be saturated at some point.
\end{remark}

Throughout the paper, we assume that agent $v_i$ detects the critical event, while the agents with subscript $j$ as $v_j\in\mathcal{V}_{\phi_i}\setminus\{v_i\}$ are asked to help agent $v_i$. Detection of a critical event by $\mathcal{D}_i^\prime$ does not necessarily mean that it is not possible to locally satisfy $\phi_i$ anymore. It rather means that the user-defined funnel boundary is touched and that repairs can help satisfying $\phi_i$. We introduce the notation $\{\hat{\boldsymbol{z}}_i\in\mathcal{Z}_i|\hat{\boldsymbol{z}}_i=\boldsymbol{z}_i\;;\;exception \}$ denoting the set of $\hat{\boldsymbol{z}}_i\in\mathcal{Z}_i$ such that $\hat{\boldsymbol{z}}_i=\boldsymbol{z}_i$ after the jump except for the elements in $\hat{\boldsymbol{z}}_i$ explicitly mentioned after the semicolon, here denoted by the placeholder $exception$.

\subsubsection{Repair of Critical Events - Stage 1}
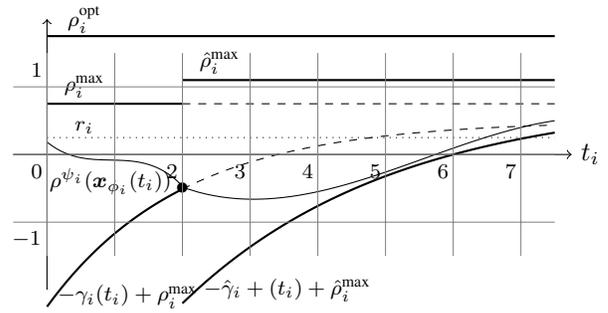
\begin{figure}
\centering
\begin{tikzpicture}[scale=0.90]
      \draw[->] (-0.5,0) -- (7.75,0) node[right] {$t_i$};
      \draw[->] (0,-2) -- (0,2) node[above] {};
      \draw[scale=1,domain=0:2,smooth,variable=\x,thick] plot ({\x},{-((3-0.25)*exp(-0.5*\x)+0.25)+0.75});
      \draw[scale=1,domain=2:7.5,smooth,variable=\x,dashed] plot ({\x},{-((3-0.25)*exp(-0.5*\x)+0.25)+0.75});
      \draw[scale=1,domain=2:7.5,smooth,variable=\x,thick] plot ({\x},{-((3.3-0.25)*exp(-0.3194*(\x-2))+0.25)+1.1});
      \draw[scale=1,domain=0:2,smooth,variable=\x,thick] plot (\x,0.75);
      \draw[scale=1,domain=2:7.5,smooth,variable=\x,dashed] plot (\x,0.75);
            \draw[scale=1,domain=-0:7.5,smooth,variable=\x,thick] plot (\x,1.75);
            \draw[scale=1,domain=2:7.5,smooth,variable=\x,thick] plot (\x,1.1);
      \draw[scale=1,domain=-0:7.5,smooth,variable=\x,dotted] plot (\x,0.25);
      \draw[scale=1,domain=-0:2,smooth,variable=\x,black] plot (\x,{-0.3*\x^3  +0.8428*\x^2  -0.8107*\x +0.184286});
            \draw[scale=1,domain=2:7.5,smooth,variable=\x,black] plot (\x,{-0.0187114*\x^3  +0.313071*\x^2  -1.38764*\x +1.18956});
      \node at (0.55,1) {\footnotesize $\rho^{\text{max}}_i$};
      \node at (0.55,0.4) {\footnotesize $r_i$};
      \node at (2.55,1.35) {\footnotesize $\hat{\rho}^{\text{max}}_i$};
      \node at (0.55,2) {\footnotesize $\rso{\psi_i}$};
      \node at (1.19,-2.1) {\footnotesize $-\gamma_i(t_i)+\rho^{\text{max}}_i$};
      \node at (3.55,-2) {\footnotesize $-\hat{\gamma}_i+(t_i)+\hat{\rho}^{\text{max}}_i$};
      \node at (0.94,-0.4) {\footnotesize $\rho^{\psi_i}(\boldsymbol{x}_{\phi_i}(t_i))$};
      \node at (-0.15,-0.25) {\footnotesize $0$};
      \node at (1.85,-0.25) {\footnotesize $2$};
      \node at (2.85,-0.25) {\footnotesize $3$};
      \node at (3.85,-0.25) {\footnotesize $4$};
      \node at (4.85,-0.25) {\footnotesize $5$};
      \node at (5.85,-0.25) {\footnotesize $6$};
      \node at (6.85,-0.25) {\footnotesize $7$};
      \node at (-0.15,1.25) {\footnotesize $1$};
      \node at (-0.3,-1.25) {\footnotesize $-1$};
      \draw (2,-0.5) node {\textbullet};
      \draw[step=1,gray,very thin] (-0.5,-1.5) grid (7.5,1.5);
\end{tikzpicture}
\caption{Funnel repair in the first stage for $\phi_i:=F_{[4,6]}\psi_i$.}
\label{fig:funnel_idea}
\end{figure}
The first repair stage is indicated by
\begin{align*}
\mathcal{D}_{i,1}^{\prime}:= &\mathcal{D}_i^\prime \cap \{(\boldsymbol{z}_i,\boldsymbol{u}^{\text{int}}_i,\boldsymbol{u}^{\text{ext}}_i)\in\mathfrak{H}_i|\mathfrak{n}_i< N_i\}
\end{align*}
where $N_i\in\mathbb{N}$ is a design parameter representing the maximum number of repair attempts in the first stage. If $(\boldsymbol{z}_i,\boldsymbol{u}^{\text{int}}_i,\boldsymbol{u}^{\text{ext}}_i)\in\mathcal{D}_{i,1}^\prime$, we first relax the parameters $t^*_i$, $\rho^{\text{max}}_i$, $r_i$, and $\gamma_i$ in a way that still guarantees local satisfaction of $\phi_i$. Pictorially speaking, we make the funnel in \eqref{eq:inequality} bigger.
\begin{example}\label{ex:repair}
Consider the formula $\phi_i:=F_{[4,6]}\psi_i$ with  $r_i:=0.25$ as the desired initial robustness, which is supposed to be achieved at $t^*_i\approx 4.8$. The original funnel is shown in  Fig. \ref{fig:funnel_idea} and given by $\rho^{\text{max}}_i$ and $-\gamma_i+\rho^{\text{max}}_i$ as in \eqref{eq:inequality}.  Without detection of a critical event, it would hence hold that $\rsct{\phi_i}{0}\ge r_i$ since $\rho^{\psi_i}\big(\boldsymbol{x}_{\phi_i}({t^*_i})\big)\ge r_i$ would be achieved.  However, at $t_r:=2$, where $t_r$ indicates the time where a critical event is detected, the trajectory $\rho^{\psi_i}\big(\boldsymbol{x}_{\phi_i}({t})\big)$ touches the lower funnel boundary and repair action is needed. This is done by setting $\hat{t}^*_i:=6$ (time relaxation), $\hat{r}_i:=0.0001$ (robustness relaxation), $\hat{\rho}^{\text{max}}_i:=1.1$ (upper funnel relaxation), and also adjusting $\hat{\gamma}_i$ (lower funnel relaxation). The funnel is hence relaxed to $\hat{\rho}^{\text{max}}_i$ and $-\hat{\gamma}_i+\hat{\rho}^{\text{max}}_i$ as depicted in Fig. \ref{fig:funnel_idea}. At the time of critical event detection $t_r$, the lower funnel is relaxed to $-\hat{\gamma}_i(t_r)+\hat{\rho}^{\text{max}}_i$ where we especially denote $\gamma^\text{r}_i:=\hat{\gamma}_i(t_r)$. Due to repair action, $\boldsymbol{x}_{\phi_i}$ locally satisfies $\phi_i$ as shown in Fig.~\ref{fig:funnel_idea}.
\end{example} 

With Example \ref{ex:repair} in mind, set
\begin{align*}
\mathcal{G}_{i,1}^\prime:=&\Big\{\hat{\boldsymbol{z}}_i\in\mathcal{Z}_i|\hat{\boldsymbol{z}}_i=\boldsymbol{z}_i\; ; \; \hat{t}^*_i := \begin{cases} b_i &\text{if } \phi_i=F_{[a_i,b_i]}\psi_i \\ 
t^*_i &\text{if } \phi_i=G_{[a_i,b_i]}\psi_i,\end{cases} \\
&\hspace{-0.8cm}\hat{\rho}^{\text{max}}_i=\rho^{\text{max}}_i+\zeta^\text{u}_i,\;  \hat{r}_i\in(0,r_i), \;\hat{\boldsymbol{p}}_i^{\gamma}=\boldsymbol{p}^{\gamma,\text{new}}_{i},\; \hat{\mathfrak{n}}_i=\mathfrak{n}_i+1\Big\}
\end{align*}
where the variables $\zeta^\text{u}_i$ and $\boldsymbol{p}^{\gamma,\text{new}}_{i}$ are defined in the sequel. In words, we set $\hat{t}_i^*:=b_i$ if $\phi_i=F_{[a_i,b_i]}\psi_i$ (time relaxation) and keep $\hat{t}_i^*:=t_i^*=a_i$ otherwise. The parameter $r_i$ is decreased to $\hat{r}_i\in (0,r_i)$ (robustness relaxation) to ensure local satisfaction of $\phi_i$. The variable $\zeta^\text{u}_i$ relaxes the upper funnel and needs to be such that $\hat{\rho}^{\text{max}}_i:=\rho^{\text{max}}_i+\zeta^\text{u}_i<\rso{\psi_i}$ (upper funnel relaxation) according to \eqref{rho_max}, i.e., let $\zeta^\text{u}_i\in(0,\rso{\psi_i}-\rho^{\text{max}}_i)$. At $t_r$, the detection time of a critical event, we set $\gamma^\text{r}_i:=\hat{\gamma}_i(t_r):=\hat{\rho}^{\text{max}}_i-\rsnt{\psi_i}+\zeta^\text{l}_i$ with 
\begin{align*}
\zeta^\text{l}_i\in\begin{cases}\mathbb{R}_{>0} &\text{if } \hat{t}^*_i>t_i\\ (0,\rsnt{\psi_i}-\hat{r}_i] &\text{otherwise}, \end{cases}
\end{align*}
which resembles \eqref{eq:g1} (lower funnel relaxation); $\zeta^\text{u}_i$ and $\zeta^\text{l}_i$ determine the margin by how much the funnel is relaxed. Let $\boldsymbol{p}^{\gamma,\text{new}}_{i}:=\begin{bmatrix} \gamma^{0,\text{new}}_{i} & \gamma^{\infty,\text{new}}_{i} & l_{i}^\text{new} \end{bmatrix}$ and select, similar to \eqref{eq:g2} and \eqref{eq:g3}, ${\gamma^{\infty,\text{new}}_{i}}\in (0,\min({\gamma^\text{r}_i},\hat{\rho}^{\text{max}}_i-\hat{r}_i)]$ and
\begin{align*}
{l_{i}^\text{new}}:=\begin{cases}
0 &\text{if } -{\gamma^\text{r}_i}+\hat{\rho}^{\text{max}}_i \ge \hat{r}_i \\
\frac{-\ln\big(\frac{\hat{r}_i+{\gamma_{i}^{\infty,\text{new}}}-\hat{\rho}^{\text{max}}_i}{-({\gamma^\text{r}_i}-{\gamma_{i}^{\infty,\text{new}}})}\big)}{\hat{t}^*_i-t_i}  &\text{if } -{\gamma^\text{r}_i}+\hat{\rho}^{\text{max}}_i< \hat{r}_i.
\end{cases}
\end{align*}
Finally, set ${\gamma^{0,\text{new}}_{i}}:=(\gamma^\text{r}_i-{\gamma^{\infty,\text{new}}_{i}})\exp({l_{i}^\text{new}}t_i)+ {\gamma^{\infty,\text{new}}_{i}}$ to account for the clock $t_i$ that is not reset $(\hat{t}_i:=t_i)$.

\subsubsection{Repair of Critical Events - Stage 2}
Repairs of the second and third stage are detected by
\begin{align*}
\mathcal{D}_{i,\{2,3\}}^{\prime}:= &\mathcal{D}_i^\prime\cap \{(\boldsymbol{z}_i,\boldsymbol{u}^{\text{int}}_i,\boldsymbol{u}^{\text{ext}}_i)\in\mathfrak{H}_i|\mathfrak{n}_i \ge N_i\}.
\end{align*}

The second stage will only be initiated if some timing constraints hold. Then, \emph{collaborative control} as in Theorem \ref{theorem:2} is used to satisfy $\phi_i$. The second stage is detected by
\begin{align*}
\mathcal{D}^\prime_{i,2}&:=\mathcal{D}_{i,\{2,3\}}^\prime\cap \Big\{(\boldsymbol{z}_i,\boldsymbol{u}^{\text{int}}_i,\boldsymbol{u}^{\text{ext}}_i)\in\mathfrak{H}_i|\forall v_j\in \mathcal{V}_{\phi_i}\setminus \{v_i\}, \\ &\hspace{-0.5cm}(\mathfrak{c}_j=-1) \text{ or } \Big(\mathfrak{c}_j=0, b_i<\begin{cases}b_j \text{ if } \phi_j=F_{[a_j,b_j]}\psi_j \\ a_j  \text{ if } \phi_j=G_{[a_j,b_j]}\phi_j \end{cases} \Big) \Big\}.
\end{align*}
i.e., each agent $v_j\in \mathcal{V}_{\phi_i}\setminus \{v_i\}$ is either free or postpones satisfaction of  $\phi_j$ to collaboratively deal with $\phi_i$ first, while ensuring that there is enough time to deal with $\phi_j$ afterwards. If $(\boldsymbol{z}_i,\boldsymbol{u}^{\text{int}}_i,\boldsymbol{u}^{\text{ext}}_i)\in\mathcal{D}^\prime_{i,2}$, all agents in $\mathcal{V}_{\phi_i}$ will use \emph{collaborative control} to deal with $\phi_i$. Therefore, let
\begin{align*}
\mathcal{G}_{i,2}^\prime:=&\Big\{\hat{\boldsymbol{z}}_i\in\mathcal{Z}_i|\hat{\boldsymbol{z}}_i=\boldsymbol{z}_i\; ; \;  \hat{\rho}^{\text{max}}_i=\rho^{\text{max}}_i+\zeta^\text{u}_i, \\
 &\hat{r}_i\in (0,r_i),\; \hat{\boldsymbol{p}}_i^{\gamma}=\boldsymbol{p}^{\gamma,\text{new}}_{i},\; \hat{\mathfrak{c}}_i=i\Big\}
\end{align*}
where $\hat{\mathfrak{c}}_i:=i$ initiates \emph{collaborative control}, while again relaxing the funnel parameters as in the first repair stage. The jump set $\mathcal{D}^\prime_{i,2}$ applies if agent $v_i$ detects a critical event. Now changing the perspective to the participating agents $v_j\in \mathcal{V}_{\phi_i}\setminus \{v_i\}$, all agents $v_j$ need to participate in \emph{collaborative control}. Assume that $v_j\in\Xi_l$, then 
\begin{align*}
\mathcal{D}^{\prime\prime}_{j,2}:=&\{(\boldsymbol{z}_j,\boldsymbol{u}^{\text{int}}_j,\boldsymbol{u}^{\text{ext}}_j)\in\mathfrak{H}_j| \mathfrak{c}_j\in\{-1,0\}, \\
&\exists v_i\in\Xi_{l}\setminus \{v_j\}, v_j\in\mathcal{V}_{\phi_i},\mathfrak{c}_i=i\},
\end{align*}
is the jump set, which is activated when agent $v_i$ asks agent $v_j$ for \emph{collaborative control}. If $(\boldsymbol{z}_i,\boldsymbol{u}^{\text{int}}_i,\boldsymbol{u}^{\text{ext}}_i)\in \mathcal{D}^{\prime\prime}_{j,2}$, set
\begin{align*}
\mathcal{G}_{j,2}^{\prime\prime}:=&\Big\{\hat{\boldsymbol{z}}_j\in\mathcal{Z}_j|\hat{\boldsymbol{z}}_j=\boldsymbol{z}_j\; ; \; \hat{\boldsymbol{p}}^\text{f}_j={\boldsymbol{p}}^\text{f}_i,\; \hat{\mathfrak{c}}_j=\mathfrak{c}_i \Big\}
\end{align*}
where $\hat{\mathfrak{c}}_j=\mathfrak{c}_i$ and $\hat{\boldsymbol{p}}^\text{f}_j={\boldsymbol{p}}^\text{f}_i$ enforce that all conditions in Theorem \ref{theorem:2} hold such that $\phi_i$ will be locally satisfied.
 
\subsubsection{Repair of Critical Events - Stage 3}

If the timing constraints in $\mathcal{D}_{i,2}^{\prime}$ do not apply, repairs of the third stage are initiated by
\begin{align*}
\mathcal{D}_{i,3}^{\prime}:= &\mathcal{D}_{i,\{2,3\}}^{\prime}\setminus \mathcal{D}_{i,2}^{\prime}.
\end{align*}

Agent $v_i$ reacts in this case by reducing the robustness $r_i$ by $\delta_i>0$ as illustrated in Example \ref{ex:1} and according to
\begin{align*}
\mathcal{G}_{i,3}^\prime:=&\Big\{\hat{\boldsymbol{z}}_i\in\mathcal{Z}_i|\hat{\boldsymbol{z}}_i=\boldsymbol{z}_i\; ; \; \hat{\rho}^{\text{max}}_i=\rho^{\text{max}}_i+\zeta^\text{u}_i, \\
 &\hat{r}_i=r_i-\delta_i,\; \hat{\rho}_i^\text{max}=\rho_i^\text{opt}+\sigma_i, \; \hat{\boldsymbol{p}}_i^{\gamma}=\boldsymbol{p}^{\gamma,\text{new}}_{i}\Big\}.
\end{align*}
where now $\gamma_i^\text{r}:=\hat{\rho}_i^\text{max}-\rho^{\psi_i}(\boldsymbol{x}_{\phi_i})+\delta_i$ is used to calculate $\boldsymbol{p}^{\gamma,\text{new}}_{i}$, while $\sigma_i>0$ will avoid Zeno behavior.

\subsubsection{The Overall System}

It now needs to be determined what happens when a task $\phi_i$ is locally satisfied. Define $\nu_i:=\begin{cases}\mathfrak{c}_i &\text{if }\mathfrak{c}_i>0  \\ i &\text{if }\mathfrak{c}_i=0 \end{cases}$ and detect such events by
\begin{small}
\begin{align*}
\mathcal{D}_{i,\text{sat}}&:=\Big\{(\boldsymbol{z}_i,\boldsymbol{u}^{\text{int}}_i,\boldsymbol{u}^{\text{ext}}_i)\in\mathfrak{H}_i|r_{\nu_i}\le \rho^{\psi_{\nu_i}}\big(\boldsymbol{x}_{\phi_{{\nu_i}}}\big)\le \rho^{\text{max}}_{\nu_i},\; \mathfrak{c}_i\ge 0,\\
&\hspace{-0.6cm}  t_i\in \begin{cases} [a_{\nu_i},b_{\nu_i}] &\text{if } \phi_{\nu_i}=F_{[a_{\nu_i},b_{\nu_i}]}\psi_{\nu_i} \\ b_{\nu_i} &\text{if } \phi_{\nu_i}=G_{[a_{\nu_i},b_{\nu_i}]}\psi_{\nu_i}\end{cases} \Big\}\setminus (\mathcal{D}^{\prime}_{i} \cup \mathcal{D}^{\prime\prime}_{i,2}),
\end{align*}
\end{small}where the set substraction of $\mathcal{D}^{\prime}_{i} \cup \mathcal{D}^{\prime\prime}_{i,2}$ exludes the case where $\mathcal{D}^{\prime}_{i}$ or $\mathcal{D}^{\prime\prime}_{i,2}$ apply simultaneously with $\mathcal{D}_{i,\text{sat}}$. This hence prevents cases when two jump options are available, which would induce an undesirable non-determism endangering the logic behind the hybrid system. If $(\boldsymbol{z}_i,\boldsymbol{u}^{\text{int}}_i,\boldsymbol{u}^{\text{ext}}_i)\in \mathcal{D}_{i,\text{sat}}$, let 
\begin{align*}
\mathcal{G}_{i,\text{sat}}:=&\Big\{\hat{\boldsymbol{z}}_i\in\mathcal{Z}_i|\hat{\boldsymbol{z}}_i=\boldsymbol{z}_i\; ; \; \hat{t}^*_i =\begin{cases} b_i &\text{if } \phi_i=F_{[a_i,b_i]}\psi_i \\ 
a_i &\text{if } \phi_i=G_{[a_i,b_i]}\psi_i,\end{cases}\\
& \hat{\rho}^{\text{max}}_i=\tilde{\rho}_i^\text{max}, \; \hat{r}_i= \tilde{r}_i,\\ & \hat{\boldsymbol{p}}_i^{\gamma}=\boldsymbol{p}^{\gamma,\text{new}}_{i}, \; \hat{\mathfrak{c}}_i=\begin{cases}0 & \text{if } \mathfrak{c}_i>0 \text{ and } \mathfrak{c}_i\neq i\\
 -1 & \text{if } \mathfrak{c}_i=0 \text{ or } \mathfrak{c}_i=i\end{cases}\Big\}
\end{align*}
where $\tilde{\rho}_i^\text{max}$ and $\tilde{r}_i$ are chosen according to \eqref{rho_max} and \eqref{r_i}, but evaluated with $\boldsymbol{x}_{\phi_i}(t_i)$ instead of $\boldsymbol{x}_{\phi_i}(0)$. If $\hat{\mathfrak{c}}_i=0$ in $\mathcal{G}_{i,\text{sat}}$, the task $\phi_i$ will be pursued next, while $\phi_i$ has already been satisfied if $\hat{\mathfrak{c}}_i=-1$ so that the agent becomes free.

Note that $\mathcal{D}_i^\prime=\mathcal{D}_{i,1}^\prime\cup\mathcal{D}_{i,2}^\prime\cup\mathcal{D}_{i,3}^\prime$ with $\mathcal{D}_{i,1}^\prime\cap\mathcal{D}_{i,2}^\prime\cap\mathcal{D}_{i,3}^\prime=\emptyset$. The hybrid system $\mathcal{H}_i$ is given by $D_i:=\mathcal{D}_i^\prime\cup\mathcal{D}_{i,2}^{\prime\prime}\cup\mathcal{D}_{i,\text{sat}}$ and $C_i := \mathcal{Z}_i\setminus D_i$. The flow map has already been defined and the jump map is
\begin{align*}
G_i&(\boldsymbol{z}_i,\boldsymbol{u}_i^{\text{int}},\boldsymbol{u}_i^{\text{ext}}):=\\
&\begin{cases}
\mathcal{G}_{i,1}^\prime(\boldsymbol{z}_i,\boldsymbol{u}_i^{\text{int}},\boldsymbol{u}_i^{\text{ext}}) &\text{for }\; (\boldsymbol{z}_i,\boldsymbol{u}_i^{\text{int}},\boldsymbol{u}_i^{\text{ext}})\in\mathcal{D}_{i,1}^\prime \\
\mathcal{G}_{i,2}^\prime(\boldsymbol{z}_i,\boldsymbol{u}_i^{\text{int}},\boldsymbol{u}_i^{\text{ext}}) &\text{for }\; (\boldsymbol{z}_i,\boldsymbol{u}_i^{\text{int}},\boldsymbol{u}_i^{\text{ext}})\in\mathcal{D}_{i,2}^\prime \\
\mathcal{G}_{i,2}^{\prime\prime}(\boldsymbol{z}_i,\boldsymbol{u}_i^{\text{int}},\boldsymbol{u}_i^{\text{ext}}) &\text{for }\; (\boldsymbol{z}_i,\boldsymbol{u}_i^{\text{int}},\boldsymbol{u}_i^{\text{ext}})\in\mathcal{D}_{i,2}^{\prime\prime} \\
\mathcal{G}_{i,3}^{\prime}(\boldsymbol{z}_i,\boldsymbol{u}_i^{\text{int}},\boldsymbol{u}_i^{\text{ext}}) &\text{for }\; (\boldsymbol{z}_i,\boldsymbol{u}_i^{\text{int}},\boldsymbol{u}_i^{\text{ext}})\in\mathcal{D}_{i,3}^{\prime} \\
\mathcal{G}_{i,\text{sat}}(\boldsymbol{z}_i,\boldsymbol{u}_i^{\text{int}},\boldsymbol{u}_i^{\text{ext}}) &\text{for }\; (\boldsymbol{z}_i,\boldsymbol{u}_i^{\text{int}},\boldsymbol{u}_i^{\text{ext}})\in\mathcal{D}_{i,\text{sat}}.
\end{cases}
\end{align*}

Note now that the sets $\mathcal{D}_{i}^\prime$ and $\mathcal{D}_{i,\text{sat}}$ as well as $\mathcal{D}_{i,2}^{\prime\prime}$ and $\mathcal{D}_{i,\text{sat}}$ are non-intersecting. However, $\mathcal{D}_{i}^{\prime}$ and $\mathcal{D}_{i,2}^{\prime\prime}$ are intersecting. Therefore, if $(\boldsymbol{z}_i,\boldsymbol{u}^{\text{int}}_i,\boldsymbol{u}^{\text{ext}}_i)\in\mathcal{D}_{i}^{\prime}\cap\mathcal{D}_{i,2}^{\prime\prime}$, which will rarely happen in practice, we only execute the jump induced by $\mathcal{D}_{i,2}^{\prime\prime}$ to not endager the logic behind the hybrid system. Thereby, we can say that the sets $\mathcal{D}_{i}^\prime$, $\mathcal{D}_{i,2}^{\prime\prime}$, and $\mathcal{D}_{i,\text{sat}}$ are technically non-intersecting.

\begin{theorem}
Assume that each agent $v_i\in\mathcal{V}$ is subject to $\phi_i$ of the form \eqref{eq:phi_class} and controlled by $\mathcal{H}_i:=(C_i,F_i,D_i,G_i)$,  while Assumptions \ref{assumption:1}-\ref{ass:5} are satisfied. The induced dependency clusters $\bar{\Xi}=\{\Xi_1,\hdots,\Xi_L\}$ are such that for each $\Xi_l\in\bar{\Xi}$ it holds that $v_i$ and $v_j$ can communicate for all $v_i,v_j\in\Xi_l$. For $v_i\in\Xi_l$ it then holds that $\rsct{\phi_i}{0}\ge r_i$, where either $r_i:=r_i(0,0)$ (initial robustness) if $\phi_i=\phi_j$ for all $v_i,v_j\in \Xi_l$  or $r_i$ is lower bounded and maximized up to a precision of $\delta_i$ otherwise. Zeno behavior is excluded.

\begin{proof}
Note first that there will never be the option of two jumps at the same time since the jump sets $\mathcal{D}_{i,1}^\prime$, $\mathcal{D}_{i,2}^\prime$, $\mathcal{D}_{i,2}^{\prime\prime}$, $\mathcal{D}_{i,3}^{\prime}$, and $\mathcal{D}_{i,\text{sat}}$ are technically non-intersecting. In the first repair stage, the parameters $t^*_i$, $\rho^{\text{max}}_i$, $r_i$, $\gamma^0_i$, $\gamma^\infty_i$, and $l_i$ are repaired in a way that still guarantees local satisfaction of $\phi_i$.  Zeno behavior is excluded for this stage since detection of a critical event is directly followed by a jump into the interior of the funnel, i.e., into the flow set $C_i$ and since only a finite number of jumps, i.e., $N_i$ jumps, are permitted. For the second repair stage, \emph{collaborative control} guarantees finishing the task $\phi_i$ by the guarantees given in Theorem \ref{theorem:2}. Afterwards, participating agents $v_j\in\mathcal{V}_{\phi_i}\setminus \{v_i\}$  have enough time to deal with their own local task $\phi_j$, which is initiated by $\mathcal{D}_{j,\text{sat}}$. If the timing constraints for \emph{collaborative control} do not hold, the third repair stage is initiated and $r_i$ is successively decreased by $\delta_i$ so that eventually $\rsct{\phi_i}{0}\ge r_i$ has to hold, i.e. maximizing $\rsct{\phi_i}{0}$ to a precision of $\delta_i$. Note that $r_i$ has to be lower bounded due to Assumption \ref{assumption:4}. This assumption states the well-posedness of $\psi_i$ and means that for local satisfaction of $\phi_i$ the state $\boldsymbol{x}_{\phi_i}$ is bounded.  Hence, all agents aim to stay within a bounded set. Consequently, successively reducing $r_i$ will eventually lead to $\rsct{\phi_i}{0}\ge r_i$. This again means that only a finite number of jumps is possible when the lower funnel is touched. Touching the upper funnel will also only lead to a finite number of jumps since $\hat{\rho}_i^\text{max}=\rho_i^\text{opt}+\sigma_i$ in $\mathcal{G}_{i,3}^\prime$, hence exluding Zeno behavior of $\mathcal{H}_i$.
\end{proof} 
\end{theorem}

\subsection{Extension to $\theta$-formulas}
If each agent $v_i\in\mathcal{V}$ is subject to $\theta_i$ of the form \eqref{eq:theta1_class}, the same result can be obtained by extending the hybrid system $\mathcal{H}_i=(C_i,F_i,D_i,G_i)$ as instructed in \cite{lindemann2017prescribed}. The detection \& repair mechanism introduced in the previous section can be applied in exactly the same way. Due to space limitations, the illustration is omitted.

\section{Simulations}
\label{sec:sim}
\begin{figure*}
\centering
\begin{subfigure}{0.48\textwidth}
\input{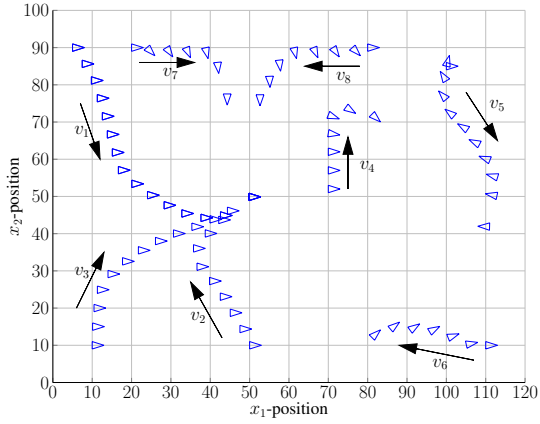}\caption{Agent trajectories for Scenario 1.}\label{fig:sim1}
\end{subfigure}
\begin{subfigure}{0.48\textwidth}
\input{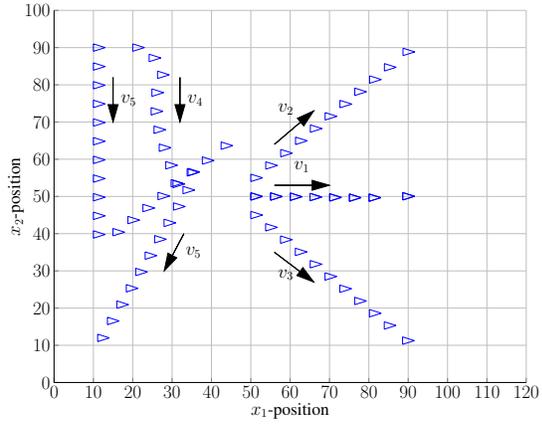}\caption{Agent trajectories for Scenario 2.}\label{fig:sim2}
\end{subfigure}
\begin{subfigure}{0.32\textwidth}
\input{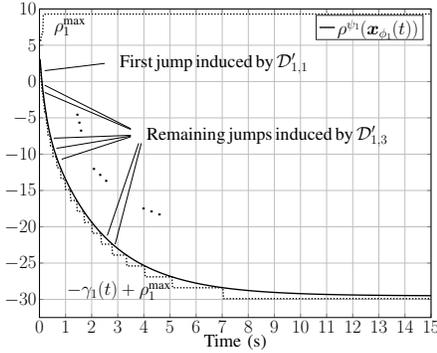}\caption{Scenario 2: Funnel repairs for agent $v_1$}\label{fig:agent1}
\end{subfigure}
\begin{subfigure}{0.32\textwidth}
\input{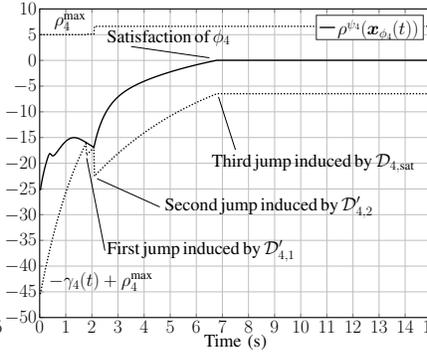}\caption{Scenario 2: Funnel repairs for agent $v_4$}\label{fig:agent2}
\end{subfigure}
\begin{subfigure}{0.32\textwidth}
\input{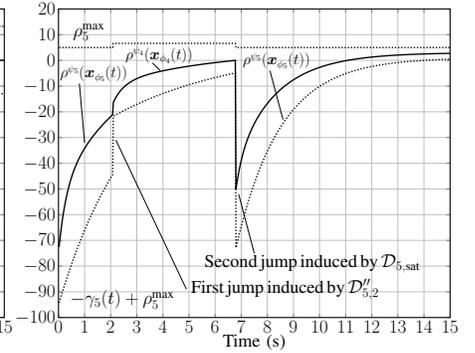}\caption{Scenario 2: Funnel repairs for agent $v_5$}\label{fig:agent3}
\end{subfigure}
\caption{Simulation results for Scenario 1 and 2.}
\end{figure*}
We consider omni-directional robots as in \cite{liu2008omni} with two states $x_1$ and $x_2$ indicating the robot position and one state $x_{3}$ indicating the robot orientation with respect to the $x_{1}$-axis. Let $x_{i,j}$ with $j\in\{1,2,3\}$ denote the $j$-th element of agent $v_i$'s state and let $\boldsymbol{p}_i:=\begin{bmatrix} x_{i,1} & x_{i,2}\end{bmatrix}$. We hence have $\boldsymbol{x}_i:=\begin{bmatrix}
\boldsymbol{p}_i & x_{i,3}
\end{bmatrix}=\begin{bmatrix}
x_{i,1} & x_{i,2} & x_{i,3}
\end{bmatrix}$ with the dynamics
\begin{align*}
\dot{\boldsymbol{x}}_i=\begin{bmatrix}
\cos(x_{i,3}) & -\sin(x_{i,3}) & 0\\
\sin(x_{i,3}) & \cos(x_{i,3}) & 0\\
0 & 0 & 1
\end{bmatrix}\Big(B_i^T\Big)^{-1}R_i\boldsymbol{u}_i,
\end{align*}
where $R_i:=0.02$ is the wheel radius and $B_i:=\begin{bmatrix}
0 & \cos(\pi/6) & -\cos(\pi/6)\\
-1 & \sin(\pi/6) & \sin(\pi/6)\\
L_i & L_i & L_i
\end{bmatrix}$ describes geometrical constraints with $L_i:=0.2$ as the radius of the robot body. Each element of $\boldsymbol{u}_i$ corresponds to the angular velocity of exactly one wheel. All simulations have been performed in real-time on a two-core 1,8 GHz CPU with 4 GB of RAM. Computational complexity is not an issue due to the computationally-efficient and easy-to-implement feedback control laws.

\emph{Scenario 1}:
This scenario employs eight agents in three clusters with $v_1,v_2,v_3\in\Xi_1$, $v_4,v_5,v_6\in\Xi_2$, and $v_7,v_8\in\Xi_3$, where the agents in each cluster are subject to the same formula, consequently  satisfying the conditions in Theorem \ref{theorem:1}. The first cluster $\Xi_1$ should eventually gather, while at the same time $v_1$ should approach the point $\boldsymbol{x}_A:=\begin{bmatrix}
50 & 50
\end{bmatrix}$. The second cluster $\Xi_2$ should eventually form a triangular formation, while the robot's orientation point to each other and agent $v_5$ approaches $\boldsymbol{x}_B:=\begin{bmatrix}
110 & 40
\end{bmatrix}$. The third cluster should eventually move to some predefined points $\boldsymbol{x}_C:=\begin{bmatrix}
40 & 70
\end{bmatrix}$ and $\boldsymbol{x}_D:=\begin{bmatrix}
55 & 70
\end{bmatrix}$, while staying as close as possible to each other and having an orientation that is pointing into the $-x_2$-direction. In formulas, we have $\phi_1:=\phi_2:=\phi_3:=F_{[10,15]}\psi_{l_1}$ with $\psi_{l_1}:=(\|\boldsymbol{p}_1-\boldsymbol{p}_2\|<2) \wedge (\|\boldsymbol{p}_1-\boldsymbol{p}_3\|<2) \wedge (\|\boldsymbol{p}_2-\boldsymbol{p}_3\|<2) \wedge (\|\boldsymbol{p}_1-\boldsymbol{p}_A\|<2)$. For the second cluster, $\phi_4:=\phi_5:=\phi_6:=F_{[10,15]}\psi_{l_2}$ is used with $\psi_{l_2}:= (\|\boldsymbol{p}_5-\boldsymbol{p}_B\|<5) \wedge (27<x_{5,1}-x_{4,1}<33) \wedge (27<x_{5,1}-x_{6,1}<33) \wedge (27<x_{4,2}-x_{5,2}<33) \wedge (27<x_{5,2}-x_{6,2}<33) \wedge (|\deg(x_{4,3})+45|<5) \wedge (|\deg(x_{5,3})-180|<5)\wedge (|\deg(x_{6,3})-45|<5)$, where $\text{deg}(\cdot)$ transforms radians into degrees. Finally, the third cluster employs $\phi_7:=\phi_8:=F_{[10,15]}\psi_{l_3}$ with $\psi_{l_3}:=(\|\boldsymbol{p}_7-\boldsymbol{p}_8\|<10) \wedge (\|\boldsymbol{p}_7-\boldsymbol{p}_C\|<10) \wedge (\|\boldsymbol{p}_8-\boldsymbol{p}_D\|<10) \wedge (|\deg(x_{7,3})+90|<5)\wedge(|\deg(x_{5,3})+90|<5)$. The simulation result is shown in Fig. \ref{fig:sim1}, where each robot's initial orientation is $0$, indicated by the direction of the triangle. Note that the tasks are satisfied within the time interval $[10,15]$.


\emph{Scenario 2}:
This scenario employs five agents in two clusters with $v_1,v_2,v_3\in\Xi_1$ and $v_4,v_5\in\Xi_2$ simulating Example \ref{ex:1} and \ref{ex:2}, respectively. Recall that in these examples we had $\phi_1:=G_{[0,15]}(\|\boldsymbol{p}_1-\boldsymbol{p}_2\|\le 10 \wedge \|\boldsymbol{p}_1-\boldsymbol{p}_3\|\le 10)$, $\phi_2:=F_{[5,15]}(\|\boldsymbol{p}_2-\begin{bmatrix}
90 & 90
\end{bmatrix}\|\le 5)$, and $\phi_3:=F_{[5,15]}(\|\boldsymbol{p}_3-\begin{bmatrix}
90 & 10
\end{bmatrix}\|\le 5)$ as well as $\phi_4:=F_{[5,10]}(\|\boldsymbol{p}_4-\boldsymbol{p}_5\|\le 10 \wedge \|\boldsymbol{p}_4-\begin{bmatrix} 50 & 70 \end{bmatrix}\|\le 10)$ and $\phi_5:=F_{[5,15]}(\|\boldsymbol{p}_5-\begin{bmatrix}
10 & 10
\end{bmatrix}\|\le 5)$ so that $\phi_1$ and $\phi_4$ are collaborative tasks. We set $\delta_i:=1.5$ and $N_i:=1$ for all agents $v_i\in\mathcal{V}$. Agent trajectories are shown in Fig. \ref{fig:sim2}, while Fig. \ref{fig:agent1} shows the funnel \eqref{eq:inequality} for agent $v_1$. It is visible that agent $v_1$ first tries to repair its parameter in Stage 1, and then initiates Stage 3 to successively reduces the robustness $r_1$ and consequently also the lower funnel as visible in Fig. \ref{fig:agent1}. Agent $v_1$ hence finds a trade-off between staying close to agent $v_2$ and $v_3$, i.e., staying in the middle of them as visible in Fig. \ref{fig:sim2}. In other words, agent $v_1$ can not satisfy $\phi_1$, but a least violating solution is found. Agent $v_4$ first tries to repair its parameters in Stage 1, but then requests agent $v_5$ to use \emph{collaborative control} to satisfy $\phi_4$ as indicated in Fig. \ref{fig:agent2} and \ref{fig:agent3}. Agent $v_5$ collaborates with agent $v_4$ to satisfy $\phi_4$ and satisfies $\phi_5$ afterwards. We can conclude that $\phi_2$, $\phi_3$, $\phi_4$, and $\phi_5$ are locally satisfied with robustness $r_2=r_3=r_4=r_5=0.5$, while $\phi_1$ is not locally satisfied, but is forced to achieve $\rsct{\phi_i}{0}>r_1=-30$.

\section{Conclusion}
\label{sec:conclusion}
We presented a framework for the control of multi-agent systems under signal temporal logic tasks. We adopted a bottom-up approach where each agent is subject to a local signal temporal logic task. By leveraging ideas from prescribed performance control, we developed a continuous feedback control law that achieves satisfaction of all local tasks under some given conditions. If these conditions do not hold, we proposed to combine the developed feedback control law with an online detection \& repair scheme, expressed as a hybrid system. This scheme detects critical events and repairs them. Advantages of our framework are low computation times and robustness that is taken care of by the robust semantics of signal temporal logic and by the prescribed performance approach.

Possible future extensions are the improvement of the repair stages in the online detection \& repair scheme. We proposed a three-stage procedure, but several other steps are possible. A next step is also to perform physical experiments.

\bibliographystyle{IEEEtran}
\bibliography{literature}

\begin{appendix}
\textit{Proof of Corollary \ref{theorem:2}:} 
We proceed similar to the proof of Theorem \ref{theorem:1}.

Step A: First, define $\boldsymbol{\xi}:=\begin{bmatrix}
\xi_{i_1} & \xi_{i_2} & \hdots & \xi_{i_{|\mathcal{V}_{\phi}|}} \end{bmatrix}$ where $|\mathcal{V}_{\phi}|$ is the cardinality of $\mathcal{V}_{\phi}$ and $v_{i_1},\hdots,v_{i_{|\mathcal{V}_{\phi}|}}\in\mathcal{V}_{\phi}$ are all agents participating in $\phi$. Define again the stacked vector $\boldsymbol{y}:=\begin{bmatrix}
\boldsymbol{x} & \boldsymbol{\xi}
\end{bmatrix}$. Consider the closed-loop system $\dot{\boldsymbol{x}}_i=:H_{\boldsymbol{x}_i}(\boldsymbol{x},\xi_i)$ with $H_{\boldsymbol{x}_i}(\boldsymbol{x},\xi_i):=f_i(\boldsymbol{x}_i)+f_i^\text{c}(\boldsymbol{x})+g_i(\boldsymbol{x}_i)\boldsymbol{u}_i+\boldsymbol{w}_i$  where
\begin{align*}
\boldsymbol{u}_i:=
\begin{cases} -\ln(-\frac{\xi_i+1}{\xi_i}){g_i(\boldsymbol{x}_i)}^T\frac{\partial \rsnt{\psi_i}}{\partial\boldsymbol{x}_i} &\text{if } v_i\in\mathcal{V}_{\phi}\\
\boldsymbol{u}_i^\prime &\text{if } v_i\in\mathcal{V}\setminus \mathcal{V}_{\phi}.
\end{cases}
\end{align*}
Recall that $\boldsymbol{u}_i^\prime$ is the control law given in the assumptions. Next, define $\dot{\boldsymbol{x}}=:H_{\boldsymbol{x}}(\boldsymbol{x},\boldsymbol{\xi})$ and $\frac{d\xi_i}{dt}=:H_{\xi_i}(\boldsymbol{x},\xi_i,t)$ with $H_{\boldsymbol{x}}(\boldsymbol{x},\boldsymbol{\xi})$ and $H_{\xi_i}(\boldsymbol{x},\xi_i,t)$ as in the proof of Theorem \ref{theorem:1}. Let
\begin{align*}
H_{\boldsymbol{\xi}}(\boldsymbol{x},\boldsymbol{\xi},t):=\begin{bmatrix} H_{\xi_{i_1}}(\boldsymbol{x},\xi_{i_1},t) & \hdots & H_{\xi_{i_{|\mathcal{V}_{\phi}|}}}(\boldsymbol{x},\xi_{i_{|\mathcal{V}_{\phi}|}},t)\end{bmatrix}
\end{align*}
so that the dynamics of $\boldsymbol{y}$ can be written as $\dot{\boldsymbol{y}}=:H(\boldsymbol{y},t)$ with
\begin{align*}
H(\boldsymbol{y},t):=\begin{bmatrix}
{H_{\boldsymbol{x}}(\boldsymbol{x},\boldsymbol{\xi})} & {H_{\boldsymbol{\xi}}(\boldsymbol{x},\boldsymbol{\xi},t)}
\end{bmatrix}.
\end{align*}  

It again holds that $\boldsymbol{x}(0)$ is such that $\xi_i(\boldsymbol{x}_{\phi_i}(0),0)\in\Omega_\xi:=(-1,0)$ holds for all agents $v_i\in\mathcal{V}_{\phi}$ by the choice of $\gamma^0_i$. As in the proof of Theorem \ref{theorem:1}, define the open, bounded, and non-empty set $\Omega_{{\phi_i}}(t)$. Next,  assume that for each agent $v_k\in\mathcal{V}\setminus \mathcal{V}_{\phi}$ the corresponding state $\boldsymbol{x}_k$ is initially contained in the open set $\Omega_{k}$, i.e., $\boldsymbol{x}_k(0)\in \Omega_{k}$,  where $\Omega_{k}$ exists due the assumptions, i.e., each state $\boldsymbol{x}_k$ remains in the  compact set $ \Omega^\prime_{k}$. Let $v_i\in\mathcal{V}_\phi$ and define 
\begin{align*}
\Omega_{\boldsymbol{\xi}}&:={\Omega_\xi}\times \hdots \times {\Omega_\xi}\subset \mathbb{R}^{|\mathcal{V}_{\phi}|}\\
\Omega_{\boldsymbol{x}}&:=\Omega_{{\phi_i}}(0)\times \Omega_{k_1} \times \hdots \Omega_{k_{M-|\mathcal{V}_{\phi}|}}\subset\mathbb{R}^{nM},
\end{align*} 
where $v_{k_1},\hdots,v_{k_{M-|\mathcal{V}_{\phi}|}}\in\mathcal{V}\setminus \mathcal{V}_{\phi}$ are all agents not belonging to $\mathcal{V}_{\phi}$. Finally, define the open, non-empty, and bounded set 
\begin{align*}
\Omega_{\boldsymbol{y}}:=\Omega_{\boldsymbol{x}} \times \Omega_{\boldsymbol{\xi}}\subset \mathbb{R}^{nM+|\mathcal{V}_{\phi}|}.
\end{align*} 
It consequently holds that $\boldsymbol{y}(0)=\begin{bmatrix}
{\boldsymbol{x}(0)} & {\boldsymbol{\xi}(0)}
\end{bmatrix}\in\Omega_{\boldsymbol{y}}$. Next, note that the conditions in Lemma \ref{theorem:sontag1} for the initial value problem $\dot{\boldsymbol{y}}=H(\boldsymbol{y},t)$ with $\boldsymbol{y}(0)\in\Omega_{\boldsymbol{y}}$ and $H(\boldsymbol{y},t):\Omega_{\boldsymbol{y}}\times \mathbb{R}_{\ge 0} \to \mathbb{R}^{nM+|\mathcal{V}_{\phi}|}$ are satisfied as in the proof of Theorem \ref{theorem:1} since the control law $\boldsymbol{u}_k^\prime$ guarantees existence of nontrivial solutions. As a result, there again exists a maximal solution with $\boldsymbol{y}(t)\in\Omega_{\boldsymbol{y}}$ for all $t\in\mathcal{J}:=[0,\tau_{\text{max}})\subseteq \mathbb{R}_{\ge 0}$ and $\tau_{\text{max}}>0$.

Step B: The Lyapunov analysis to show that $\tau_{\text{max}}=\infty$ follows similar steps as in the proof of Theorem \ref{theorem:1} that are not shown here. It can again be shown that $\xi_i(t)\in\Omega_{\xi_i}^\prime$ and $\boldsymbol{x}_{\phi_{i}}(t)\in \Omega_{{\phi_{i}}}^\prime$ for all agents $v_i\in\mathcal{V}_{\phi}$, where 	$\Omega_{\xi_i}^\prime$ and 	$\Omega_{{\phi_{i}}}^\prime$ are compact subsets of $\Omega_{\xi_i}$ and 	$\Omega_{{\phi_{i}}}$, respectively. It consequently follows that $0<r\le\rsct{\phi}{0}\le\rho^{\text{max}}$, i.e., $(\boldsymbol{x}_\phi,0)\models \phi$. 
\end{appendix}

\addtolength{\textheight}{-12cm}   

\end{document}